\documentclass{amsart}
\usepackage{amscd,amsfonts,amsmath,amssymb,amsthm,eucal,latexsym} 
\usepackage[dvips]{graphicx}

\def\nd{\noindent}

\def\iT#1{{\itshape#1\/}}

\numberwithin{equation}{section}
\newtheorem*{thm*}{Theorem} 
\newtheorem{thm}{Theorem}
\newtheorem{prop}{Proposition}
\newtheorem{lem}[prop]{Lemma}
\newtheorem{cor}[prop]{Corollary}
\theoremstyle{definition}
\newtheorem{definition}[prop]{Definition} 
\newtheorem{ex}[prop]{Example}

\newtheorem{rem}[prop]{Remark}
\newtheorem{conjecture}[prop]{Conjecture}
\newtheorem{notation}[prop]{Notation}
\numberwithin{prop}{section}

\def\Rem#1{Remark~\ref{#1}}
\def\Def#1{Definition~\ref{#1}}
\def\Cor#1{Corollary~\ref{#1}}
\def\Lem#1{Lemma~\ref{#1}} 
\def\Thm#1{Theorem~\ref{#1}} 
\def\Sec#1{Section~\ref{#1}} 
\def\Fig#1{Figure~\ref{#1}} 
\def\Prop#1{Proposition~\ref{#1}} 
 
\def\Conj#1{Conjecture~\ref{#1}} 
\def\Ex#1{Example~\ref{#1}} 
 
\def\Sec#1{Section~\ref{#1}}

\long\def\BTHMm#1#2{\begin{thm*}[#1]#2\end{thm*}}
\long\def\BTHM#1#2{\begin{thm}\LL{#1}#2\end{thm}}
\long\def\BDF#1#2{\begin{definition}\LL{#1}#2\end{definition}}
\long\def\BPROP#1#2{\begin{prop}\LL{#1}#2\end{prop}}
\long\def\BLEM#1#2{\begin{lem}\LL{#1}#2\end{lem}}

\long\def\BREM#1#2{\begin{rem}\LL{#1}#2\end{rem}}
\long\def\BEX#1#2#3{\begin{ex}[#2]\LL{#1}#3\end{ex}}
\long\def\BCOR#1#2{\begin{cor}\LL{#1}#2\end{cor}}
\long\def\BNOT#1#2{\begin{notation}\LL{#1}#2\end{notation}}
\long\def\BCONJ#1#2{\begin{conjecture}\LL{#1}#2\end{conjecture}}

\long\def\BEQ#1#2{\begin{equation}\LL{#1}#2\end{equation}}
 
\long\def\BFIG#1#2#3{\begin{figure}[th]#3\caption{#2}\LL{#1}\end{figure}} 

\def\Pr{\begin{proof}}
\def\rP{\end{proof}}

\def\boxx#1{\leavevmode\vbox{\hbox to 0pt{\hss\raise1.8ex\vbox 
to 0pt{\vss\hrule\hbox{\vrule\kern.75pt\vbox{\kern.75pt\hbox{\tiny #1}\kern.75pt}\kern.75pt\vrule}\hrule}}}\relax} 

\def\LL#1{\label{#1}\protect\boxx{#1}} 
\def\LL#1{\label{#1}}

\def\iFF{{if and only if }} 
\def\ph{{piecewise harmonic}} 
\def\pcf{{p.c.f.}} 
\def\pcfsss{{\pcf\ self-similar set}} 
\def\s-s{{self-similar}} 
\def\fr{finitely ramified} 
\def\Fr{Finitely ramified} 
\def\fr{fi\-ni\-te\-ly ram\-i\-fied} 
\def\frf{\fr\ fra\-ct\-al} 
 
\def\frcs{\fr\ cell struct\-ure}

\def\Str{Strichartz}

\def\Sig{Sier\-pi\'n\-ski gasket}

\def\Sic{Sierpi\'nski carpet}
\def\Lp{Laplacian}
\def\Df{Dirichlet form}
\def\hc{harmonic coordinates}

\def\a{\alpha}
\def\A{\mathcal A}

\def\E{\mathcal E}
\def\F{\mathcal F}

\def\L{\mathcal L}

\def\V{\mathcal V}
\def\x{(x)}
\def\<{\langle}
\def\>{\rangle}
\def\Tr{\text{\rm Tr\,}}
\def\Trace{\text{\rm Trace}\vph_}
\def\Grad{\text{\rm Grad\,}}

\def\Tan{\text{\rm Tan\,}}
\def\Tana{\text{\rm Tan}_\alpha}
\def\Dom{\text{\rm Dom\,}}
\def\vph{^{\vphantom{{\ }_{Ap}}}}

\def\hhh{\put(0,0){\mbox{see  \texttt{\SMALL http://www.math.uconn.edu/$\sim$teplyaev/research/arxiv0506261.pdf}}}}

\begin{document}

\title
[Harmonic coordinates on fractals]{Harmonic coordinates on fractals \\ with finitely ramified cell structure}
\author{Alexander Teplyaev}
\thanks{Research supported in part by the NSF grant DMS-0505622}

\address{\noindent Department of Mathematics, 
University of Connecticut, Storrs CT 06269-3009 USA}
\email{teplyaev@math.uconn.edu}
\date{\today}
\begin{abstract}  
We define sets with finitely ramified  cell structure, which are generalizations of \mbox{p.c.f.} self-similar 
sets introduced by Kigami and of fractafolds introduced by Strichartz. In general, 
we do not assume even local self-similarity, and allow countably many cells 
connected at each junction point. 
In particular, we consider post-critically infinite fractals. 
We prove that if Kigami's resistance form 
satisfies certain assumptions, then there exists a weak Riemannian metric 
such that the energy can be expressed as the integral of the norm squared 
of a weak gradient with respect to an energy measure. 
Furthermore, we prove that if such a set can be homeomorphically represented 
in harmonic coordinates, then for smooth functions the weak gradient can be 
replaced by the usual gradient. 
We also prove a simple formula for the energy measure Laplacian in harmonic 
coordinates. 
{\tableofcontents}
\end{abstract}


\maketitle

\section*{Introduction} 
There is a well developed theory 
of Dirichlet (energy, resistance) forms, and corresponding random processes, on the class of 
post-critically finite (\pcf\ for short) \s-s\ sets, which are 
 finitely ramified \cite[and references therein]{Ba,Ki1,Ki,St2,St4}. Also, many piecewise and stochastically self-similar fractals have been considered \cite[and references therein]{BH,H,H2,St3}. The general non \s-s\ energy forms 
on the \Sig\ were studied in \cite{MST}. 
In all the mentioned works the fractals considered have \frcs. In this paper we will 
extend some aspects of this theory for a  class of spaces, which may have no self-similarity in any sense, 
and may have infinitely many cells   connected at every junction point. 
Throughout this paper we extensively and substantially use the general theory of resistance 
forms developed in \cite{Ki4}. The 
existence 
of such forms is a delicate question even in the self-similar 
\pcf\ case \cite[and references therein]{HMT,Ki4,Me}. 
To prove our results we use some methods introduced in \cite{T1}. 
In the present paper we give the basic background information, and the reader may find all the details in \cite{Ki4,T1}. 

 In 
\Sec{section-Resistance-K} we give the definition of a resistance form in the sense of Kigami \cite{Ki4}. 
In \Sec{section-Definitions} we define sets with \frcs s. 
Examples of such fractals are 
 \pcf\ self-similar sets introduced by Kigami in \cite{Ki1,Ki}, fractafolds introduced by \Str\ in \cite{St3.5}, random fractals 
\cite[and references therein]{BH,H,H2}, and 
 non self-similar \Sig s \cite{MST,T2}. The key topological assumption is that there 
is a cell structure such that every cell has finite boundary, but we do not assume any self-similarity. 

The terminology we use can be explained as follows. 
The term ``post-critically infinite'' in this context means that every junction point can be an intersection of countably infinite number of cells with pairwise disjoint interior, that is every cell can be linked to 
 countably many other cells. 
The term  ``\fr'' means that every cell is joined with its complement in a finite number of points. 
A good example of an infinitely ramified fractal 
 is the \Sic. There exists a \s-s\ diffusion and corresponding 
Dirichlet form on the \Sic\ \cite{BB1,BB2,BB3,KZ}, but its uniqueness has not been proved. 

In 
\Sec{section-Resistance} we prove that  Kigami's resistance form  is a local regular Dirichlet form under appropriate conditions. 
In \Sec{section-Riemannian}
 we prove that if the resistance form satisfies certain non degeneracy assumptions, then there exists a weak Riemannian metric, defined almost everywhere, such that the energy can be expressed as the integral of the norm of a weak gradient with respect to an energy measure. This generalizes earlier results by Kusuoka \cite{Ku2} and the author \cite{T1}. 
Furthermore, in \Sec{section-Gradient} we prove that if the \frf\ can be homeomorphically represented in harmonic coordinates, then the weak gradient can be replaced by the usual gradient for smooth functions, which generalizes an earlier result by Kigami in \cite{Ki2}. 
In \Sec{section-Lp} we prove a simple formula for the energy measure Laplacian in harmonic coordinates. 
This formula was announced, in the case of the standard energy 
form on the \Sig, in \cite{T2} without a proof. In a sense, the formula for the energy \Lp\ is the second derivative with respect to the generalized Riemannian metric. 
In the case of the standard energy 
form on the \Sig, it is proved by Kusuoka in \cite{Ku1} that this generalized Riemannian metric has rank one almost everywhere. This can be interpreted as that in harmonic coordinates on the \Sig\ the energy \Lp\ is the one dimensional second derivative in the tangential direction. 
We conjecture that this is the case for any \frf\ 
considered in this paper. The main tool we use in this theorem is  approximating the \frf\ by a sequence of so called quantum graphs \cite[and references therein]{Kuchment1,Kuchment2}. 
In \Sec{section-self-similar} we discuss \s-s\ \frf s, and existence of \s-s\ resistance forms in particular. 
In \Sec{section-Examples} we give several examples of \frf s for which our theory can be applied. Among them are factor-spaces of \pcfsss s, and  post-critically infinite analogs of the \Sig.

In the case of the standard energy 
form on the \Sig, it is proved by Kigami in \cite{Ki7} that the heat kernel with respect to the energy measure has Gaussian asymptotics in harmonic coordinates (a weaker version was obtained in \cite{MS}). Recently a powerful machinery was developed to obtain heat kernel estimates on various ``rough'' spaces, including many fractals \cite[and references therein]{BBK,Ki6}. It is not unlikely that this theory is applicable to many, if not all, \frf s in \hc. 
Also, 
 some results about the singularity of the 
energy measure with respect to product measures \cite{BST,Hino1,Hino2} are valid in the case of \fr\ \s-s\ fractals under suitable extra assumptions.

\subsubsection*{Acknowledgments} 
The author is very grateful to Richard Bass, Klara Hveberg, Jun Kigami and Robert \Str\ for 
many important remarks and suggestions.

\section{Kigami's resistance forms}\label{section-Resistance-K}


 Below we restate the definition of a resistance form in \cite{Ki4}. 

\BDF{d-resistance form}{ A pair $(\E,\Dom\E)$ is called a resistance form on a countable set $V_*$ if it satisfies the following conditions.
\begin{itemize}
\item[(RF1)] $\Dom\E$ is a linear subspace of $\ell(V_*)$ containing constants, 
 $\E$ is a nonnegative symmetric quadratic form on $\Dom\E$, and $\E(u,u)=0$ if and only if $u$ is constant on $V_*$.
\item[(RF2)] Let $\sim$ be the equivalence relation on $\Dom\E$ defined by $u\sim v$ \iFF $u-v$ is constant on $V_*$. Then $(\E/\mbox{\hskip-.4em}\sim,\Dom\E)$ is a Hilbert space. 
\item[(RF3)] For any finite subset $V\subset V_*$ and for any 
$v\in\ell(V)$ there exists $u\in\Dom\E$ such that $u\big|_V=v$.
\item[(RF4)] For any $p,q\in V_*$ 
$$
\sup\left\{ \frac{\big(u(p)-u(q)\big)^2}{\E(u,u)}:u\in\Dom\E,\E(u,u)>0\right\}<\infty.
$$
This supremum is denoted by $R(p,q)$ and called the effective resistance between $p$ and $q$. 
\item[(RF5)] For any $u\in\Dom\E$ we have the $\E(\bar u,\bar u)\leqslant\E(u,u)$, where 
$$
\bar u(p)=\left\{
\begin{aligned}
&1&& \text{if \,} u(p)\geqslant1,\\
&u(p)&& \text{if \,} 0<u(p)<1,\\
&0&& \text{if \,} u(p)\leqslant1.
\end{aligned}\right.
$$
Property~(RF5) is called the Markov property. 
\end{itemize}}

Note that the effective resistance $R$ is a metric on $V_*$, 
and that any function in $\Dom\E$ is $R$-continuous. Let $\Omega$ 
be the $R$-completion of $V_*$. Then  any $u\in\Dom\E$ 
has a unique $R$-continuous extension to $\Omega$. 

For any finite subset $U\subset V_*$ the 
finite dimensional Dirichlet form $\E_U$ on $U$ is defined by 
$$\E_U(f,f)=\inf\{\E(g,g):g\in\Dom\E, g\big|_U=f\},$$ 
which exists by  \cite{Ki4}, and moreover there is a unique $g$ for which the $\inf$ is attained. 
The Dirichlet form $\E_U$ is called the trace of $\E$ on $U$, and denoted $$\E_U=\Trace U(\E).$$ By the definition, if $U_1\subset U_2$ then 
 $\E_{U_1}$ is  the trace of $\E_{U_2}$ on $U_1$, that is $\E_{U_1}=\Trace {U_1}(\E_{U_2})$. 

\BTHMm{Kigami \cite{Ki4}}{Suppose that $V_n$ are finite subsets of $V_*$ and that $\bigcup_{n=0}^\infty V_n$ is $R$-dense in $V_*$. Then  
  $$\E(f,f)=\lim_{n\to\infty}\E_{{V_n}}(f,f)$$ for any $f\in\Dom\E$, where the limit is actually non-decreasing. Is particular, $\E$ is uniquely defined by the sequence of its finite dimensional traces $\E_{V_n}$ on ${V_n}$.}

\BTHMm{Kigami \cite{Ki4}}{Suppose that  $V_n$ are finite sets,     for each $n$ there is a resistance form $\E_{V_n}$ on ${V_n}$, and this sequence of finite dimensional forms is compatible in the sense that 
each $\E_{V_n}$ is  the trace of $\E_{V_{n+1}}$ on ${V_n}$, where $n=0,1,2,...$. Then there exists a resistance form $\E$ on $V_*=\bigcup_{n=0}^\infty V_n$ such that 
  $$\E(f,f)=\lim_{n\to\infty}\E_{{V_n}}(f,f)$$ for any $f\in\Dom\E$, and the limit is actually non-decreasing.}

\section{\Fr\ fractals} \label{section-Definitions}

\BDF{d-fract}{A \iT{\frf} $F$ is 
a compact metric space with a \iT{cell structure}~$\F=\{F_\alpha\}_{\alpha\in\A}$ and 
a \iT{boundary (\text{vertex}) structure} $\V=\{V_\alpha\}_{\alpha\in\A}$ 
such that 
the following conditions hold.
\begin{itemize}
\item[(A)] 
$\A$ is a countable index set; 
 \item[(B)] 
each $F_\alpha$ 
is a distinct compact connected subset of $F$; 
 \item[(C)] 
 each  $V_\alpha$ is a finite subset of 
$F_\alpha$ with at least two elements;
 \item[(D)] 
if $F_\alpha=\bigcup_{j=1}^k F_{\alpha_j}$  then  
$V_\alpha\subset\bigcup_{j=1}^k V_{\alpha_j}$; 
 \item[(E)] 
{there exists a filtration $\{\A_n\}_{n=0}^\infty$ such that 
\begin{enumerate}
\item 
$\A_n$ are finite subsets of $\A$,  $\A_{0}=\{0\}$, and 
$F_{0}= F$; 
\item 
$\A_n\cap\A_m=\varnothing$ if $n\neq m$; 
 \item 
for any $\alpha\in\A_n$ there are  
${\alpha_1},...,{\alpha_k}\in\A_{n+1}$
such that $F_\alpha=\bigcup_{j=1}^k F_{\alpha_j}$; 
\end{enumerate}}
 \item[(F)] 
$F_{\alpha'}\bigcap F_{\alpha\vph}=V_{\alpha'}\bigcap V_{\alpha\vph}$ for any 
two distinct $\alpha,\alpha'\in\A_n$; 
\item[(G)] 
for any strictly decreasing infinite cell sequence $F_{\alpha_1}\supsetneq F_{\alpha_2}\supsetneq...$ there exists $x\in F$ such that 
 $\bigcap_{n\geqslant1} F_{\alpha_n} =\{x\}$. 
\end{itemize}
If these conditions are satisfied, then $$(F,\F,\V)=
(F,\{F_\alpha\}_{\alpha\in\A},\{V_\alpha\}_{\alpha\in\A})$$ 
is called a {\it \frcs.}}

\BNOT{not-Vn}{We denote $V_n=\bigcup_{\alpha\in\A_n}V_\alpha $. 
Note that $V_n\subset V_{n+1}$ for all $n\geqslant0$ by \Def{d-fract}.
We say that $F_\alpha$ is an $n$-cell if $\alpha\in \A_n$.}

\BREM{rem-pcf}{By the definition every cell in a \frf\ has a  boundary consisting of isolated points (see \Prop{prop-fundamental}), which implies the name ``point connected''. In particular, any \pcfsss\ is a \frf. However, every vertex $v\in V_*$ of a \frf\ can be 
 an intersection of countably many cells with pairwise disjoint interior (see \Ex{e-nSG}). Hence even if a \frf\ is self-similar, it does not have to be a \pcfsss.}

\BREM{rem-V0}{In this definition the vertex boundary 
$V_{0}$ of $F_{0}=F$ can be arbitrary, and in general may 
have no relation with the topological structure of $F$. 
However, the (WN) and (HC) assumptions made below will 
de facto impose restrictions on the choice of $V_{0}$. In particular, 
the energy measure, gradient and the energy measure \Lp\ all depend on the choice of~$V_{0}$. 

This is somewhat different from the theory of \pcfsss s in \cite{Ki1,Ki3,Ki}, where 
$V_{0}$ is uniquely determined as the post-critical set of the \pcf\ \s-s\ structure. Note, however, that the same topological fractal $F$ can have different \s-s\ structures, and different post-critical sets in particular.

Note that every \pcfsss\ is a \frf, see \cite{Ki1,Ki} and \Sec{section-self-similar}.}

\BREM{rem-filtration}{In general a filtration is not unique for a \frf. 
For example, the filtration $\widetilde\A_{k}=\A_{2k}$ satisfies all the conditions of \Def{d-fract}.
However, the results of this paper do not depend on  the choice of the filtration. In particular, Kigami's resistance forms, energy
measures \mbox{etc.} are independent of the filtration. If the
fractal is self-similar, then changing the self-similar structure sometimes is
very useful, and changing filtration is one of the results of changing it. Moreover, if the \frf\ is not
self-similar, then it may not be clear what the natural filtration is.
\emph{In what follows we assume that some 
filtration is fixed.}}

\BPROP{prop-x}{For any $x\in F$ there is a strictly decreasing infinite sequence of cells satisfying condition (G) of the definition. The diameter of cells in any such sequence tend to zero. } 

\Pr Suppose $x\in F$ is given. We choose $F_{\alpha_1}=F$. Then, if $F_{\alpha_n}$ is chosen, we choose $F_{\alpha_{n{+}1}}$ to be a proper sub-cell of $F_{\alpha_n}$ which contains~$x$.  
Suppose for a moment that the diameter of cells in  such a sequence does not tend to zero. 
Then  for each $n$ there is $x_n\in F_{\alpha_n}$ such that $\liminf_{n\to\infty} d(x_n,x)=\varepsilon>0$. By compactness there is $y\in \bigcap_{n\geqslant1} F_{\alpha_n}$ such that $d(y,x)\geqslant\varepsilon$. 
 This is a contradiction with the property~(G) of \Def{d-fract}. \rP

\BPROP{prop-fundamental-1}{The topological boundary of $F_\alpha$ 
is contained in $V_\alpha$ for any $\alpha\in\A$.} 

\Pr For any closed set $A$ we have $\partial A=A\cap Closure(A^c)$, where $A^c$ is the complement of $A$. If $A=F_\alpha$ is an $n$-cell, then $Closure(A^c)$ is the union of all $n$-cells except $F_\alpha$. Then the proof 
 follows from property~(F) of \Def{d-fract}. \rP

\BPROP{prop-count}{The set $V_*=\bigcup_{\alpha\in\A} V_\alpha$ is countably infinite, and $F$ is uncountable.} 

\Pr The set $V_*$ is a countable union of finite sets, and every cell is a union of at least two smaller sub-cells. Then each cell is uncountable by properties (B) and (C) of \Def{d-fract}.\rP

\BPROP{prop-fundamental-0}{For any distinct $x,y\in F$ there is $n(x,y)$ such that if $m\geqslant n(x,y)$ 
then any $m$-cell can not contain both $x$ and $y$. } 

\Pr Let $\mathcal B_m(x,y)$ be the collection of all $m$-cells that contain both $x$ and $y$.
By  definition any cell in $\mathcal B_{m+1}(x,y)$ is contained in a cell which  belongs to $\mathcal B_m(x,y)$. Therefore,   if there are infinitely many nonempty collections $\mathcal B_m(x,y)$, then there is an infinite  decreasing sequence of cells that contains 
 both $x$ and $y$. This is a contradiction with  property~(G) of \Def{d-fract}. \rP

\BPROP{prop-fundamental}{For any $x\in F$ and $n\geqslant 0$, let $U_n(x)$ denote the union of all $n$-cells that contain $x$. Then the collection of open sets 
$\mathcal U=\{U_n(x) ^\circ\}_{x\in F,n\geqslant 0}$ is a countable fundamental sequence of neighborhoods. Here $B^\circ$ denotes the topological interior of a set $B$. 

Moreover, for any $x\in F$ and open neighborhood $U$ of $x$ there exist $y\in V_*$ and $n$ 
such that $x\in U_n(x)\subset U_n(y)\subset U$. In particular,  the smaller collection of open sets 
$\mathcal U'=\{U_n(x) ^\circ\}_{x\in V_*,n\geqslant 0}$ is  a countable fundamental sequence of neighborhoods.} 

\Pr Note that the collection $\mathcal U'$  is countable because $V_*$ is countable by \Prop{prop-count}. 
The collection $\mathcal U$ is countable  because if $x$ and $y$ belong to the interior of the same $n$-cell, 
then $U_n(x)=U_n(y)$. 

First, suppose $x\in V_*$. Then we have to show that   for any  open neighborhood $U$ of $x$ there exists 
 $n\geqslant 0$ such that $ U_n(x)\subset U$. 
Suppose for a moment that such $n$ does not exist. Then for any $n$ the set $U_n(x)\backslash U$ 
is a nonempty compact set. Moreover, the sequence of sets $\{U_n(x)\backslash U\}_{n\geqslant 0}$ is decreasing 
and so has a nonempty intersection. Then we can choose $z\in \bigcap_{n\geqslant 0} U_n(x)\backslash U  $, and for any $n$ there is an $n$-cell that contains both $x$ and $z$. This is a contradiction with \Prop{prop-fundamental-0}.

Now suppose $x\notin V_*$. Then  for any $n>0$ 
 there exists $y_n\in V_n$ such that 
$x\in U_n(y_n)\subset U_{n-1}(x)$. Moreover, we can assume also that $U_n(y_n)\subset U_{n-1}(y_{n-1})$  for any $n>1$. 
Then we have to show that   for any  open neighborhood $U$ of $x$ there exist 
 $n>0$ such that $ U_n(y_n)\subset U$. 
Suppose for a moment that such $n$ does not exist. Then the set $U_n(y_n)\backslash U$ 
is a nonempty compact set. Moreover, the sequence of sets $\{U_n(y_n)\backslash U\}_{n\geqslant 1}$ is decreasing 
and so has a nonempty intersection. Then we can choose $z\in \bigcap_{n\geqslant 1} U_n(y_n)\backslash U  $, and for any $n>1$ there is an $(n-1)$-cell that contains both $x$ and $z$. This is a contradiction with \Prop{prop-fundamental-0}. \rP

\section{Resistance forms on \frf s}\label{section-Resistance}

We assume that there is a resistance form  on $V_*$ in the sense of Kigami \cite{Ki,Ki4}, see 
\Def{d-resistance form}.
For convenience we will denote $\E_n(f,f)=\E_{V_n}(f,f)$, see \Sec{section-Resistance-K}.  Recall that 
$\E(f,f)=\lim_{n\to\infty}\E_{n}(f,f)$ for any $f\in\Dom\E$, 
where the limit is actually non-decreasing. 

\BDF{def-harm}{A function is harmonic if it minimizes the energy for the given set of boundary values.}
Note that any harmonic function is uniquely defined by its restriction to $V_0$. Moreover, any 
function on  $V_0$ has a unique continuation to a harmonic function. For any harmonic function $h$ we have $\E(h,h)=\E_n(h,h)$ for all $n$ by \cite{Ki4}.
Also note that for any function $g\in\Dom\E$ we have $\E_0(g,g)\leqslant\E(g,g)$, and a function $h$ is harmonic if and only if  $\E_0(h,h)=\E(h,h)$.

Let $\E_\alpha(f,f)=(\E_\alpha)_{V_\alpha}(f,f)$, where $\E_\alpha$ is the restriction of $\E$ to 
$F_\alpha$. 
Then  $$\E_n=\sum_{\alpha\in\A_n}\E_{V_\alpha}.$$

\BLEM{prop-atoms}{If $h$ is harmonic and  continuous then 
$$\lim_{n\to\infty}\sum_{\alpha\in\A_n,x\in F_\alpha}\E_{\alpha}(h\big|_{V_\alpha},h\big|_{V_\alpha})=0  $$
for any $x\in F$.}

\Pr Let  $\E(h,h)=e>0$. It is easy to see that the limit under consideration is decreasing and so it 
 exists. 
Suppose for a moment this limit is equal to $c>0$. 

Without loss of generality we can assume that $h\x=0$ and that $|h(y)|\geqslant1$ for any $y\in V_0\backslash\{x\}$. 
By \Prop{prop-x} for 
any $\varepsilon>0$ there are cells $F_{\a_1},...,F_{\a_l}$ such that $|h(x)-h(y)|<\varepsilon$ for any $y\in\bigcup_{j=1}^lF_{\a_j}$, and $\bigcup_{j=1}^lF_{\a_j}$ 
contains a neighborhood of $x$. Without loss of generality we 
can assume that $V_0\bigcap\left(\bigcup_{j=1}^lF_{\a_j}\backslash\{x\}\right)=\varnothing$. 

Let $V'=\bigcup_{j=1}^lV_{\a_j}$ and consider the trace of the resistance 
form on $V_0\bigcup V'$. Obviously if $\varepsilon$ is small then there is 
a uniform bound for conductances between points in $V_0\backslash\{x\}$ 
and $ V'$. Then consider changing the values of $h$ on $V'$ to zero. 
Inside of $\bigcup_{j=1}^lF_{\a_j}$ the energy will be reduced by at 
least $c$, since the function is now constant there. On the other hand, 
outside of $\bigcup_{j=1}^lF_{\a_j}$ the energy increase will be bounded 
by a constant times $\varepsilon e$. So the total energy will decrease if 
$\varepsilon$ is small enough. This is a contradiction with the 
definition of a harmonic function, and so $c=0$. 

Note that the proof works even if $V'$ is an infinite set and so 
it is applicable to connected spaces with cell structure, such as the \Sic, which is not a \frf.\rP 

\BCOR{cor-atoms}{If $h$ is harmonic and  continuous then there is a unique continuous energy measure 
$\nu_h$ on $F$ defined by $\nu_h(F_\alpha)=\E_{\alpha}(h\big|_{V_\alpha},h\big|_{V_\alpha})$ for all $\alpha\in\A$.}

\BREM{rem-h-cont0}{\emph{In what follows we assume that harmonic functions are continuous.}}

\BDF{d-nu}{We fix a complete, up to constant functions, energy orthonormal set of harmonic functions $h_1,...,h_k$, 
where $k=|V_0|-1$, and define the Kusuoka energy measure by 
$$
\nu=\nu_{h_1}+...+\nu_{h_k}.
$$}

If $F_{\alpha'}\subset F_{\alpha\vph}$, then 
$$
M_{\alpha\vph{,}\alpha'}:
\ell{(V_{\alpha\vph})}\to\ell{(V_{\alpha'})}
$$ 
is the linear map which is defined as follows. If $f_\a$ is a function on $V_\a$ then let $h_{f_\a}$ be the unique harmonic function on $F_\a$ that coincides with $f_\a$ on $V_\a$. Then we define 
$$
M_{\alpha\vph{,}\alpha'}f_\a=h_{f_\a}\big|_{V_{\a'}}.
$$ 
Thus $
M_{\alpha\vph{,}\alpha'}
$
transforms the (vertex) boundary values of a harmonic function on $F_{\alpha\vph}$ 
into the values of this harmonic function on $V_{\alpha'}$. 
We denote $
M_{\alpha}\vph=M_{0,\alpha}\vph$. 
We denote $D_\alpha$ the matrix of the  Dirichlet form $\E_\alpha$ on $V_\alpha$. 
By elementary linear algebra  we have the following lemma  (see \cite{T1} and 
also \cite{Ki1,Ki,Ku1}). 

\BLEM{prop-prop}{If $F_\alpha=\bigcup F_{\alpha_j}$ then 
$$
D_\alpha\vph=
\sum
M_{\alpha{,}\alpha_j}^*D_{\alpha_j}\vph M_{\alpha{,}\alpha_j}\vph$$
and 
$$
\nu(F_\alpha\vph)=
\Tr
M_{\alpha}^*D_{\alpha}\vph M_{\alpha}\vph
.$$
In particular $\nu$ is defined uniquely in the sense that it does not depend on the choice of the complete energy orthonormal set of harmonic functions.}

We denote $$\displaystyle
Z_\alpha\vph=
\frac{M_{\alpha}^*D_{\alpha}\vph M_{\alpha}\vph}
{\nu(F_\alpha\vph)}$$ 
if $\nu(F_\alpha\vph)\neq0$. 
Then we define matrix valued functions 
$$Z_n\x=Z_\alpha$$ 
if $\nu(F_\alpha\vph)\neq0$, $\alpha\in\A_n$ and $x\in F_\alpha\backslash V_\alpha$. Note that $\Tr Z_n\x=1$ by definition. 

\BTHM{thm-Z}{For $\nu$-almost all $x$ there is a limit $$Z\x=\lim_{n\to\infty} Z_n\vph\x.$$}

\Pr 
One can see, 
following the original Kusuoka's idea \cite{Ku1,Ku2}, 
that $Z_n$ is a bounded $\nu$-martingale.
\rP

\BREM{rem-Z-dnu}{In a sense, the matrix valued measure  $Zd\nu$ plays the role of a generalized Riemaninan metric on the fractal $F$ (see Theorems~\ref{thm-C1} and~\ref{thm-C2}). The matrix $Z$ has trace one by its definition, but on many fractals it is discontinuous.  Moreover, in some  examples, such as the \Sig, the matrix $Z$ has rank one almost everywhere. Then it can be described as the projection onto the one dimensional tangent space.

One can see that the energy measures $\nu_h$ are the same as the energy measures in the 
 general theory of \Df s \cite{BouH,FOT}. 
One can also define the matrix $Z$ as the matrix whose entries are the densities 
$$Z_{ij}=\frac{d\nu_{h_i,h_j}}{d\nu}$$
using the general theory of Dirichlet forms in \cite{BouH,FOT}. However we give a different description because the pointwise approximation using the cell structure is important in this theorem. }

\BDF{def-n-harm}{A function is  $n$-harmonic if it minimizes the energy for the given set of  values on $V_n$.} 
Note that any $n$-harmonic function is uniquely defined by its restriction to $V_n$. Moreover, any 
function on  $V_n$ has a unique continuation to an $n$-harmonic function. Also note that for any function $g\in\Dom\E$ we have $\E_n(g,g)\leqslant\E(g,g)$, and a function $f$ is $n$-harmonic if and only if  $\E_n(f,f)=\E(f,f)$.

Recall that $R$ is the effective resistance metric on $V_*$, 
and that any function in $\Dom\E$ is $R$-continuous. Let $\Omega$ 
be the $R$-completion of $V_*$. Then  any $u\in\Dom\E$ 
has a unique $R$-continuous extension to $\Omega$. The next theorem generalizes \cite[Proposition 3.3.2]{Ki} for 
possibly non \s-s \frf s.

\BTHM{thm-theta}{Suppose that all $n$-harmonic functions are continuous. Then any continuous function is $R$-continuous, and  any $R$-Cauchy sequence converges in the topology of $F$. Also, there is a continuous 
injective map $\theta:\Omega\to F$ which is the 
identity on $V_*$.}

\Pr It is easy to see from the maximum principle that any continuous function can be uniformly approximated by $n$-harmonic functions, which implies that any continuous function is $R$-continuous. Suppose for a moment that $\{x_k\}$ is an $R$-Cauchy sequence in $V_*$ which does not converge.  By compactness, it must have a limit point, say $x$. Then, by  the results of \Sec{section-Definitions}, there is $n$ and two disjoint $n$-cells $F_{\alpha}$ and $F_{\beta}$ such that 
$x\in F_{\beta}$, but $F_{\alpha}$ contains an infinite subsequence of $\{x_k\}$, say $\{y_m\}$. There is an $n$-harmonic function $f$ which is identically 1 on $F_{\beta}$ and zero on every point of $V_n$ which is not in $F_{\beta}$. Then for any $m$ we have $R(x,y_m)\geqslant1/\E(f,f)>0$, which is a contradiction. Thus,  any $R$-Cauchy sequence converges in the topology of $F$. Therefore we can define a  continuous 
 map $\theta:\Omega\to F$ which is the 
identity on $V_*$. Now Suppose for a moment  that $\theta$ is not injective. Then there are two 
$R$-Cauchy sequences in $V_*$, say $\{x_k\}$  and $\{y_k\}$, which  have the same limit in $F$ but two different $R$-limits in $\Omega$, say $x$ and $y$. By continuity, for any $m$-harmonic function $f$ we have $f(x)=f(y)$. This is a contradiction since the space of $m$-harmonic functions separates points of $\Omega$ by \cite{Ki4}. Thus, $\theta$ is  injective. 
\rP 

\BREM{rem-h-cont}{If conditions of \Thm{thm-theta} are satisfied, then we can (and will) consider 
$\Omega$ as a subset of $F$. Then $\Omega$ is the 
$R$-closure of $V_*$. In a sense, $\Omega$ is the set where 
 the Dirichlet form $\E$ ``lives".}

\BTHM{thm-locreg}{Suppose that all $n$-harmonic functions are continuous. Then  $\E$ is a local regular \Df\ on $\Omega$ (with respect to any measure that charges every nonempty open set).}

\Pr The regularity of $\E$ is proved in \cite{Ki4}. In particular, $\Dom\E\mod(constants)$ is a Hilbert space in the energy norm. Note that the set of $n$-harmonic functions is a core of $\E$ in both the original  and $R$-topologies. Also note that if a set is $R$-compact then it is compact in the original topology of $f$ by \Thm{thm-theta}. 
Suppose now $f$ and $g$ are two functions in $\Dom\E$ with disjoint compact supports. Then, by the results of \Sec{section-Definitions}, there is $n$ and a finite number of $n$-cells $F_{\alpha_1},...,F_{\alpha_k}$ such that $\bigcup_{i=1}^kF_{\alpha_i}$ contains the support of $f$ but is disjoint with the support of $g$. Then it is easy to see that for any $m\geqslant n$ we have $\E_m(f,g)=0$ and so $\E(f,g)=0$.
\rP

\section{Generalized Riemannian metric and weak gradient}\label{section-Riemannian}

\BDF{def-loc-harm}{We say that $f\in \Dom\E$ is $n$-\ph\ if 
 for any $\alpha\in\A_n$ there is a (globally) harmonic 
function $h_\alpha$ that coincides with $f$ on $F_\alpha$.}

Note that, by definition, the notion of $n$-\ph\ functions in general is more 
restrictive than the more commonly used notion of $n$-harmonic 
functions  defined in the previous section.

\BDF{def-weakly-nondegenerate}{We say that the resistance form on a \frf\ is \emph{weakly nondegenerate} if the space of \ph\ functions is dense in $\Dom\E$.} 

The notion of weakly nondegenerate harmonic structures was studied in \cite{T1} in the case of \pcfsss s. 

\medskip 

\nd{\bf Assumption (WN).} \emph{In what follows we assume that the resistance form is weakly nondegenerate.} 

\BPROP{prop-WN-nu}{The (WN) assumption implies $\text{supp}(\nu)=F$.}

\Pr Our definitions imply that for any cell $F_\a$ there is a 
function of finite energy with support in this cell. If it can be approximated by \ph\ functions, then $\nu(F_\a)>0$. \rP

By \Prop{prop-fundamental} $\text{supp}(\nu)=F$ if and only if every cell has a positive measure. 

\BTHM{prop-supp-nu}{Let $F_\nu$ be the factor-space (quotient) of $F$ obtained by collapsing all cells of zero $\nu$-measure. Then $F_\nu$ 
 is a \frf\ with the cell and vertex structures naturally inherited from $F$. }

\Pr The only nontrivial condition to verify is that any cell of $F_\nu$ has at least two boundary points. 
The maximum principle implies that a cell $F_\alpha$ has a positive $\nu$-measure 
if and only if there is a harmonic function which is non constant on $V_\alpha$. \rP

\BDF{def-Tan}{If $f$ is $n$-\ph\ then we define 
its tangent $\Tana f$ for $\alpha\in\A_n$ as the unique element of $\ell(V_0)$ 
that satisfies two conditions: 
\begin{itemize}
\item[(A)] 
if $h_{\alpha,\Tan}$ is the harmonic function with boundary values $\Tana f$ then $h_{\alpha,\Tan}$ coincides with $f$ on $F_\alpha$; 
 \item[(B)] 
$h_{\alpha,\Tan}$ has the smallest energy among all harmonic functions $h_\alpha$ such that $h_{\alpha}$ coincides with $f$ on $F_\alpha$. 
\end{itemize}}

We define $L^2_Z$ as the Hilbert space of $\ell(V_0)$-valued 
functions on $F$ with the norm defined by 
\[
\|u\|_{L^2_Z}^2=\int_ F 
\< u , Z u \> d\nu . 
\]

\BDF{def-Grad}{If $f$ is $n$-\ph\ then we define 
its gradient $\Grad f$ as the element of $L^2_Z$ such that, for $\nu$-almost all $x$, 
$\Grad f\x=\Tana f$ in the sense of $L^2_Z$ if $x\in F_\alpha$ and $\alpha\in\A_n$.}

\BLEM{lem-grad}{If $f$ is $n$-\ph\ then $\E(f,f)=\|\Grad f\|_{L^2_Z}^2$.}

\Pr Follows from \Lem{prop-prop}.\rP

\BTHM{thm-Grad}{Under the (WN) assumption $\Grad$ can be extended 
from the space of \ph\ functions to an isometry $$Grad: \Dom\E\to\L^2_Z,$$ which is called the weak gradient. }

\Pr The statement follows from \Lem{lem-grad} and the (WN) assumption. \rP

\BCOR{thm-nu-f}{Under the (WN) assumption we have 
$$\nu_ f\ll\nu$$ for any $f\in \Dom\E$.} 
 
\Pr The statement follows from \Thm{thm-Grad}. It can also be obtained directly from the (WN) assumption, or the general theory of \Df s \cite{BouH,FOT}.\rP

\BCONJ{conj-EM-WN}{We conjecture that the assumption $\text{supp}(\nu)=F$ is equivalent to the (WN) assumption for all 
\frf s.}

\BCONJ{conj-Z}{We conjecture that for any  \frf\ $\text{rank}Z\x=1$ for $\nu$-almost all $x$. }

The next proposition follows easily from our definitions. It means, in particular, that \Conj{conj-Z} implies \Conj{conj-EM-WN}.

\BPROP{prop-conj}{If $\text{supp}(\nu)=F$ and 
$\text{rank}Z\x=1$ for $\nu$-almost all $x$ then the (WN) assumption holds.}

\section{Gradient in \hc}\label{section-Gradient}

To define harmonic coordinates one needs to choose a complete, up to constant functions, set of harmonic functions $h_1,...,h_k$ and define the coordinate map $\psi:F\to{} \mathbb R^k$ by 
$\psi(x){=}(h_1\x,... ,h_k\x )$. 
A particular choice of harmonic coordinates is not 
important since they are equivalent up to a linear change of variables. 
Below we fix the most standard coordinates which make the computations simpler. 

\BDF{def-HC}{Let $V_0=\{v_1,...,v_m\}$ and let $h_j$ 
be the unique harmonic function with boundary values 
$h_j(v_i)=\delta_{i,j}$. 
Kigami's harmonic coordinate map 
$\psi:F\to{} \mathbb R^m$ is defined by 
$\psi(x){=}(h_1\x,... ,h_m\x )$.} 

\BLEM{prop-HC-}{\begin{enumerate}
\item 
Any set $\psi(F_\alpha)$ is contained in the convex hull of $\psi(V_\alpha)$. 
\item 
A set $\psi(F_\alpha)$ has at least two points if and only if $\psi(V_\alpha)$ 
has at least two points. 
\item \label{i-i-i}
If on 
$F_ H=\psi(F)$  we define a cell structure that consists of all sets $\psi(F_\alpha)$ that have at least two points, then conditions (A)--(E) and (G) of \Def{d-fract} are satisfied.
\item 
If for all $n$ and  for any 
two distinct $\alpha,\alpha'\in\A_n$ we have 
$$\psi(F_{\alpha'})\bigcap \psi(F_{\alpha\vph})=\psi(V_{\alpha'})\bigcap \psi(V_{\alpha\vph}),$$ then 
$F_ H=\psi(F)$  is a \frf\ with the cell structure  defined in Item (\ref{i-i-i}) of this lemma. 
\end{enumerate}}

\Pr The maximum principle implies that $\psi(F_\alpha)$ is contained in the convex hull of $\psi(V_\alpha)$, which implies the other statements. \rP

The next theorem easily follows from this lemma. 

\BTHM{cor-HC}{ $\psi:F\to F_ H=\psi(F)$ is a homeomorphism if and only if 
for any $\alpha\in\A$ the map $\psi\big|_{V_\alpha}$ 
is an injection, and $$\psi(F_{\alpha'}\cap F_{\alpha\vph})=
\psi(F_{\alpha'})\cap \psi(F_{\alpha\vph})$$for all $\alpha,\alpha'\in\A$. }

\nd{\bf Assumption (HC).} \emph{In what follows we assume that 
$\psi:F\to F_ H=\psi(F)$ is a homeomorphism.}

\BPROP{prop-HC}{The (HC) assumption implies the (WN) assumption.} 

\Pr It is easy to see that under the (HC) assumption any cell has 
 positive measure, and that any continuous function can be uniformly approximated by \ph\ functions. The latter is true because all harmonic functions are linear in harmonic coordinates, and the maximum principle implies that $\psi(F_\alpha)$ is contained in the convex hull of $\psi(V_\alpha)$. \rP 

\BNOT{not-ell}{In what follows, for simplicity, we assume $F= F_H$ 
and $\psi\x=x $. Also, we identify $\ell(V_0)$ with $\mathbb R^m$ 
in the natural way.}

\BTHM{thm-C1}{Under the (HC) assumption we have that if 
$f$ is the restriction to $F$ of a 
$C^1(\mathbb R^m)$ function then $f\in \Dom\E$, and such functions are dense 
in $\Dom\E$. Moreover, if $f\in C^1(\mathbb R^m)$ then
$$\Grad f=\nabla f$$ 
in the sense of the Hilbert space $L^2_Z$. 
In particular we have 
 the 
 \iT{Kigami formula} 
\[
\E(f,f)=\|\nabla f\|_{L^2_Z}^2= \int_ F 
\< \nabla f , Z \nabla f \> d\nu 
\]for any $f\in C^1(\mathbb R^m)$. }

\Pr In fact, we will prove this result for a somewhat larger space of functions. We say that $f$ is a piecewise $C^1$-function if for some $n$ and for all $\alpha\in\A_n$ there is $f_\alpha\in C^1(\mathbb R^m)$ such that 
$f_\alpha \big|_{F_\alpha}=f\big|_{F_\alpha}$. 
In particular, a \ph\ function is piecewise $C^1$. 

If $g$ is a linear function in $R^m$ then 
$g \big|_{V_0}=\nabla g$ 
 since we identify 
$\ell(V_0)$ with $\mathbb R^m$ in the natural way. 
Therefore for 
any \ph\ function $f$ we have $\Grad f=\nabla f$ in the sense of the Hilbert space $L^2_Z$. 

Any $C^1$-function is a  
piecewise $C^1$-function, and any piecewise $C^1$-function can be approximated by \ph\ (that is, piecewise linear) functions in 
$C^1$ norm. 
Thus, to complete the proof we need an estimate of the energy of a function in terms of its $C^1$ norm, provided by the next simple \Lem{lem-C1}. \rP

\BLEM{lem-C1}{If $f$ is the restriction to $F$ of a 
$C^1(\mathbb R^m)$ function then 
\BEQ{e-est}{\E_n(f,f)\leqslant\nu(F)\|f\|_{C^1(\mathbb R^m)}^2}and the same estimate holds for $|\E(f,f)|$.}

\Pr By definition \cite{Ki1,Ki} of $\E_n$  we have that 
\BEQ{e-est-}{\begin{aligned}&\E_n(f,f)=\sum_{x,y\in V_n}c_{n,x,y}\big(f(x)-f(y)\big)^2\leqslant\\&
\|f\|_{C^1(\mathbb R^m)}^2\sum_{x,y\in V_n}c_{n,x,y}|x-y|^2=
\|f\|_{C^1(\mathbb R^m)}^2\nu(F).\end{aligned}}\rP

\BREM{rem-C1}{Using \Thm{thm-Grad} one can prove 
\Thm{thm-C1} using the general theory of Dirichlet forms in \cite{BouH,FOT} (see \Rem{rem-Z-dnu}).
However we give a constructive proof which also 
defines an approximating sequence to the gradient. A similar proof can be made along the lines of the proof of \Thm{thm-C2} using approximations by quantum graphs. }

\section{Energy measure \Lp\ in harmonic coordinates}\label{section-Lp}

The the energy measure \Lp\ can be defined as follows. We say that  \mbox{$f\in\Dom\Delta_\nu$} if there exists a function  $\Delta_\nu f\in L^2_\nu $ such that  
\BEQ{eq-Gauss-Green}{\E(f,g)=-\int_Fg\Delta_\nu f d\nu,} 
for any function  $g\in\Dom\E$  vanishing on the boundary $V_0$. By \cite{Ki4}, the \Lp\ $\Delta_\nu$ 
is a uniquely defined linear operator with    $\Dom\Delta_\nu\subset\Dom\E$. In fact 
$\Dom\Delta_\nu$ is $\E$-dense in $\Dom\E$, and is also dense in $L^2_\nu $.  The \Lp\ $\Delta_\nu$ 
is  self-adjoint with, say,  Dirichlet or Neumann boundary conditions. Formula \eqref{eq-Gauss-Green} 
is often called  the Gauss-Green formula. Extensive information on the relation of a \Df\ and its generator, the \Lp, can be found in \cite{BouH,FOT,Ki}.

\BTHM{thm-C2}{Under the (HC) assumption we have that if 
$f$ is the restriction to $F$ of a 
$C^2(\mathbb R^m)$ function then 
$f\in \Dom{\Delta_\nu }$, and such functions are $\E$-dense 
in $\Dom{\Delta_\nu }$. Moreover, $\nu$-almost everywhere 
$$\Delta_\nu f = \Tr(Z D^2f )$$ 
where $D^2f$ is the matrix of the second derivatives of $f$.}

\def\QQ{^{\mbox{\,\hskip-.5ex}^{Q}}}
\Pr 
We start with defining a different sequence of approximating energy forms. In various situations these forms are associated with so called quantum graphs, photonic crystals and cable systems. 
If $f\in C^1(\mathbb R^m)$ then we define 
\[
\E_n\QQ(f,g)=\sum_{x,y\in V_n}c_{n,x,y}\E_{x,y}\QQ(f,f)
\]
where 
\[
\E_{x,y}\QQ(f,f)=\int_0^1\Big(\tfrac{d}{dt}f\big(x(1-t)+ty\big)\Big)^2dt
\]
is the integral of the square of the derivative 
$$
\tfrac{d}{dt}f\big(x(1-t)+ty\big)=\<\nabla f\big(x(1-t)+ty\big),y-x\>
$$ 
of $f$ along the straight line 
segment connecting $x$ and $y$. Thus 
$\E_{x,y}\QQ(f,f)$ is the usual one dimensional 
energy of a function on a straight line segment. If $f$ is linear then 
$\E_{x,y}\QQ(f,f)=\big(f(x)-f(y)\big)^2$. Therefore if $f$ is \ph\ then 
$ \E_n\QQ(f,f)=
\E_n(f,f)$ for all large enough $n$. Also $\E_n\QQ$ satisfies estimate \eqref{e-est}. Therefore for any $C^1(R^m)$-function we have 
$$
\lim_{n\to\infty}\E_n\QQ(f,f)=\E(f,f)
$$
by \Thm{thm-C1}. 

It is easy to see that if $g$ is a $C^1(R^m)$-function vanishing on $V_0$ and 
$f$ is a $C^2(R^m)$-function then 
\[
\E_n\QQ(f,g)=\sum_{x,y\in V_n}c_{n,x,y}
\int_0^1
g\big(x(1-t)+ty\big)\Big(\tfrac{d^2}{dt^2} f\big(x(1-t)+ty\big)\Big) dt
\]
because after integration by parts all the boundary terms are canceled. 
Then if $\alpha\in\A_n$ then 
\[\begin{aligned}&
\sum_{x,y\in V_\alpha}c_{n,x,y}
\tfrac{d^2}{dt^2} f\big(x(1-t)+ty\big) =
\\&
\sum_{x,y\in V_\alpha}c_{n,x,y}\sum_{i,j=1}^m
{D^2_{ij}} f\big(x(1-t)+ty\big)(y_i-x_i)(y_j-x_j) =
\\&
\Tr\Big(M_\alpha^*D\vph_\alpha M\vph_\alpha\big(D^2 f(x_\alpha) + R_n(x,y,t,f,\alpha,x_\a)\big)\Big)
\end{aligned}\]
where $x_\alpha\in V_\alpha$ and 
$$
\lim_{n\to\infty}|R_n(x,y,t,f,\alpha,x_\a)|=0
$$
uniformly in $\a\in\A_n$, $x,y,x_\a\in F_\a$ and $t\in[0,1]$, 
which completes the proof.\ Note also that one can obtain an estimate similar to \eqref{e-est}, as in \Cor{lem-C2}.\rP

\BCOR{thm-C2-0}{Under the (HC) assumption, 
$\Delta_\nu f \in L^\infty(F)$ for any 
$f\in C^2(\mathbb R^m)$.}

\BCOR{cor-C2-1}{Under the (HC) assumption, if 
$f(x)=\|x\|^2$ then 
$\Delta_\nu f = 1$.}

\BCOR{lem-C2}{If $f$ is the restriction to $F$ of a 
$C^2(\mathbb R^m)$ function, and $g$ is the restriction to $F$ of a 
$C^1(\mathbb R^m)$ function vanishing on the boundary, then 
$$|\E_n(f,g)|\leqslant\nu(F) \|g\|\vph_{C(\mathbb R^m)}\|f\|\vph_{C^2(\mathbb R^m)}$$
and the same estimate holds for $|\E(f,g)|$.}

\Pr This estimate follows from the proof of \Thm{thm-C2}.\rP

\BREM{rem-C2}{One can also obtain \Thm{thm-C2} from \Thm{thm-C1} using the general theory of Dirichlet forms in \cite{BouH,FOT} (see \Rem{rem-Z-dnu}). However we give a different constructive proof using the approximation by 
quantum graphs (see \cite{Kuchment1,Kuchment2}). }

\section{Topologically \s-s \frf s}\label{section-self-similar}

\BDF{definition-sspcnf}{A compact connected metric space  $F$ is called a \iT{\fr\ \s-s set} if 
there are injective contraction maps $$\psi_1,...,\psi_m:F\to F$$ and a finite set $V_0\subset F$  
such that $$F=\bigcup_{i=1}^m\psi_i(F)$$ and for any $n$ and for any two distinct words $w,w'\in W_n=\{1,...,m\}^n$ we have 
$$F_w\cap F_{w'}=V_w\cap V_{w'},$$ 
where $F_w=\psi_w(F)$ and $V_w=\psi_{w}(V_0)$. 
Here for a finite word $w=w_1...w_n\in W_n$ we denote 
 $$\psi_w=\psi_{w_1}\circ...\circ\psi_{w_n}.$$   
The set $V_0$ is called the vertex boundary of $F$. }

\BPROP{prop-sspcnf}{A \fr\ \s-s set is a \frf\ provided $V_0$ has  at least two elements. 

We have 
$\A_n=W_n$ for $n\geqslant1$ and $\A=\{0\}\bigcup W_*$, where $W_*=\bigcup_{n\geqslant1 }W_n$.}

\Pr All items  in  \Def{d-fract}  are self-evident. Note that item (B) holds because each cell is connected and has at least two elements, and the intersection of two cells is finite. 
Item (G) holds because $\psi_i$ are contractions. 
\rP

\BREM{rem-Hveberg}{The question of existence of a ``self-similar'' metric 
on self-similar sets was recently studied in detail in \cite{Hveberg}. 
According to \cite{Hveberg}, our class of \s-s \frf s defined 
above is the same as finitely ramified SSH-fractals (with finite fractal boundary) of \cite{Hveberg}. 
The definition of SSH-fractals in \cite{Hveberg} requires fulfillment of a certain set of axioms, 
one of which is that the maps $\psi_1,...,\psi_m:F\to F$  
are continuous injections. It is then proved that $F$ can be equipped with 
a self-similar metric in such a way that the
injective maps  $\psi_j$  become contractions (as well as local
similitudes), but the topology does not change. 
We use a simplified approach when we assume from the beginning that $\psi_i$ are contractions. 

In addition, it is proved in \cite{Hveberg} that for every \pcfsss\ defined in  \cite{Ki1,Ki} 
there exists a  self-similar metric.  Therefore our definition of \s-s \frf s generalizes the definition of \pcfsss s. 
Our definition allows 
 infinitely many cells to meet at a junction point, which is  referred to as fractals with ``infinite
multiplicity'' in  \cite{Hveberg}.}

Note that, by definition, each $\psi_i$ maps $V_*$ into itself injectively. 

\BDF{d-SSE}{A resistance form $\E$ on $V_*$, in the sense of \Sec{section-Resistance}, is \s-s with energy  renormalization factors $\rho=(\rho_1,...,\rho_m)$  if for any $f\in\Dom\E$ we have 
\BEQ{e-rho}{\E(f,f)=\sum_{i=1}^m\rho_i\E(f_i,f_i).}
Here we use the notation $f_w=f\circ\psi_w$ for any $w\in W_*$. }

The energy renormalization factors, or weights, $\rho=(\rho_1,...,\rho_m)$ are often also called conductance 
scaling factors because of the relation of resistance forms and electrical networks. They are reciprocals of the 
resistance scaling factors $r_j=\frac1{\rho_j}$.

\BDF{def-Lambda}{For a set of energy  renormalization factors $\rho=(\rho_1,...,\rho_m)$ and any resistance form $\E_0$ on $V_0$ define the resistance form $\Psi_\rho(\E_0)$ on $V_1$ by 
$$\Psi_\rho(\E_0)(f,f)=\sum_{i=1}^m\rho_i\E_0(g_i,g_i),$$
where
$$g_i=f\big|\vph_{\psi_i(V_0)}\circ\psi^{-1}_i.$$
Then 
$\Lambda(\E_0)$ is defined as the trace of $\Psi_\rho(\E_0)$  on $V_0$:
$$\Lambda(\E_0)=\Trace{V_0}\Psi_\rho(\E_0).$$}

The next two propositions are essentially proved in \cite{Ki,Ki4,Me}.

\BPROP{prop-Lambda1}{If $\E$ is \s-s then $\E_0=\Lambda(\E_0)$.}

\BPROP{prop-Lambda2}{If $\E_0$ is such that $\E_0=\Lambda(\E_0)$ then there is a \s-s resistance form $\E$ such that $\E_0$ is the trace of $\E$ on $V_0$.}

\BTHM{thm-hf-cont}{On any \s-s\ \frf\ with a \s-s resistance form all $n$-harmonic functions are continuous.}

\Pr By the \s-sity, it is enough to prove that the harmonic functions are continuous. Since all $\psi_i$ 
are contractions, there is $n$ such that any $n$-cell contains at most one point of $V_0$. 
By the strong maximum principle there is $\varepsilon>0$ such that for any $w\in W_n$ and any harmonic function $h$ we have $$\left|\max_{x\in F}h(x)-\min_{x\in F}h(x)\right|\geqslant(1-\varepsilon)\left|\max_{x\in F_w}h(x)-\min_{x\in F_w}h(x)\right|.$$ Then for any positive integer $m$ and any $w\in W_{mn}$  we have $$\left|\max_{x\in F}h(x)-\min_{x\in F}h(x)\right|\geqslant(1-\varepsilon)^m\left|\max_{x\in F_w}h(x)-\min_{x\in F_w}h(x)\right|.$$ 
\rP

We conjecture that the results of \cite[Section 3.3]{Ki}, and many other results of \cite{Ki,Ki4} on the topology and analysis on \pcfsss\ hold for \fr\ \s-s sets as well. The next theorem is one of these results. 
Following \cite{Ki1,Ki}, we say that the \s-s resistance form is \emph{regular}\/ if $\rho_i>1$ for all $i$. 

\BTHM{thm-reg}{If a \s-s resistance form on a \s-s\ \frf\ $F$ is regular, then $\Omega=F$.}

\Pr If $\text{diam}\vph_R(\cdot)$ denotes the diameter of a set in the effective resistance metric $R$, and $\rho_w=\rho_{w_1}...\rho_{w_n}$ 
for any finite word $w=w_1...w_n\in W_n$, then $$\text{diam}\vph_R(F) \geqslant \rho_w\text{diam}\vph_R(F_w)$$ 
by the self-similarity of the resistance form and the definition of the metric $R$.
\rP

\BDF{def-group}{The group $G$ is said to act on a \frf\ $F$ if each $g\in G$ is a homeomorphism of $F$ such that  $g(V_n)=V_n$ for all $n\geqslant0$. }

\BPROP{prop-group}{If a group $G$  acts on a \frf\ $F$ then for each $g\in G$ and each $n$-cell $F_\alpha$, $g(F_\alpha)$ is an $n$-cell.
}

\Pr From the results of \Sec{section-Definitions} we have  that $n$-cells are connected, have pairwise disjoint interiors, and their topological boundaries are contained in $V_n$, which is preserved by $g$ by definition. 
\rP

\BTHM{thm-symm}{Suppose a group $G$  acts on a \s-s \frf\ $F$ and $G$ restricted to $V_0$ is the whole permutation group of $V_0$. Then there exists a unique, up to a constant,  $G$-invariant  \s-s resistance form $\E$ with equal energy  renormalization weights and \BEQ{e-e0}{\E_0(f,f)=\sum_{x,y\in V_0}\big(f(x)-f(y)\big)^2.}}

\Pr 
It is easy to see that, up to a constant, $E_0$ is the only $G$-invariant resistance form on $V_0$. Let $\rho_1=(1,...,1)$. Then $\Lambda(\E_0)$ is also $G$-invariant and so $\E_0=c\Trace{V_0}\Psi_{\rho_1}(\E_0)$ for some $c$. Then the result holds for $\rho=c\rho_1$ by \Prop{prop-Lambda1} and \Prop{prop-Lambda2}. 
\rP

An $n$-cell  is called  a  boundary cell if it intersects
$V_0$. Otherwise it is called an interior cell. 
We say that $F$ has  connected interior if
the set of interior $1$-cells is connected, any boundary $1$-cell contains exactly one point of $V_0$, and
the intersection of two different boundary $1$-cells is  contained in
an interior $1$-cell.
The following theorem is proved in \cite{HMT} for the \pcf\ case, but the proof applies for \s-s \frf\ without any changes. 

\BTHMm{\cite{HMT}}{Suppose that $F$ has connected interior, and a group $G$  acts on a \s-s \frf\ $F$ such that its action on $V_0$ is 
transitive. Then there exists a  $G$-invariant  \s-s resistance form $\E$.}

Other results in \cite{HMT} also apply for \s-s \frf.

\section{Examples}\label{section-Examples}

\BEX{e-I}{Unit interval}{The usual unit interval is a \frf. In this case $V_*$ can be any countable dense subset of [0,1] which includes $\{0,1\}$. The usual energy 
form $$\E(f,f)=\int_0^1|f'(t)|^2dt$$ satisfies all the assumptions of our paper. The energy measure is the 
Lebesgue measure and the \Lp\ is the usual second derivative. }

\BEX{e-Q}{Quantum graphs}{A quantum graph, a collection of finite number of points in $\mathbb R^m$ joined by weighted straight line segments (see \cite{Kuchment1,Kuchment2} and also the proof of \Thm{thm-C2}), is 
a \frf. The usual energy 
form on a quantum graph, which is the sum of weighted standard one dimensional forms on each segment, satisfies all the assumptions of our paper.}


\def\LiNe#1#2#3#4#5{{
\count41=#3 \advance\count41 by #1 
\count42=#4 \advance\count42 by #2 
\divide\count41 by 2 
\divide\count42 by 2 
\qbezier(#1,#2)(\count41,\count42)(#3,#4)}}



\def\LiNe#1#2#3#4#5{{
\count101=#1 
\count102=#2 
\count41=#3 \advance\count41 by-\count101 
\count42=#4 \advance\count42 by-\count102 
\count51=\count41 \ifnum\count41<0 \multiply\count41 by-1 \fi
\count52=\count42 \ifnum\count42<0 \multiply\count42 by-1 \fi 
\advance\count41 by \count42 \divide\count41 by \dvv 
\advance\count41 by 1
\divide\count51 by \count41 \divide\count52 by \count41 
\advance\count41 by 1 
\multiput(\count101,\count102)(\count51,\count52){\count41}{#5}
}}


\def\LiNe#1#2#3#4#5{{
\count41=#3 \advance\count41 by-#1 
\count42=#4 \advance\count42 by-#2 
\count51=\count41 \ifnum\count41<0 \multiply\count41 by-1 \fi
\count52=\count42 \ifnum\count42<0 \multiply\count42 by-1 \fi 
\advance\count41 by \count42 \divide\count41 by \dvv 
\advance\count41 by 1 
\count42=\count41 
\count43=0 
\loop 
\count101=#1 
\count102=#2 
\count103=#3 
\count104=#4 
\multiply\count101 by \count42 \multiply\count103 by \count43 
\advance\count101 by \count103
\divide\count101 by \count41 
\multiply\count102 by \count42 \multiply\count104 by \count43 
\advance\count102 by \count104
\divide\count102 by \count41 
\put(\count101,\count102){#5}
\advance\count42 by-1 
\advance\count43 by 1 
\ifnum \count42>-1 \repeat
}}



\def\trigg#1#2#3#4#5#6#7{{\LiNe{#1}{#2}{#3}{#4}{#7} 
\LiNe{#3}{#4}{#5}{#6}{#7}\LiNe{#5}{#6}{#1}{#2}{#7}}}

\def\trigg#1#2#3#4#5#6#7{{
\count107=#1 
\count108=#2 
\advance\count107 by #3 \advance\count107 by #5 \divide\count107 by 3 
\advance\count108 by #4 \advance\count108 by #6 \divide\count108 by 3
\LiNe{#1}{#2}{\count107}{\count108}{#7} 
\LiNe{#3}{#4}{\count107}{\count108}{#7} 
\LiNe{#5}{#6}{\count107}{\count108}{#7} 
}}


\def\power#1#2{\count91=1 \count92=#2 \loop 
\advance\count92 by-1 \multiply\count91 by #1 
\ifnum \count92>0 \repeat}

\def\troo#1{\divide\count101 by 2 \divide\count102 by 2 
\ifcase#1 
 \or 
\advance\count101 by 8192 \or 
\advance\count101 by 4096 \advance\count102 by 7168 \fi} 

\def\SIGA#1#2{ \power{3}{#1}
\loop 
\advance\count91 by-1 
\count92=#1 \count95=\count91 
\count101=0 \count102=0 
{\loop 
\advance\count92 by-1 
\count93=\count95 
\divide\count95 by 3 
\count94=\count95 
\multiply\count94 by 3 
\advance\count93 by-\count94 
\troo{\count93} 
\ifnum \count92>0 \repeat \put(\count101,\count102){#2}}
\ifnum \count91>0 \repeat} 

\def\trooh#1{
\ifcase#1 
%
%
%
\count141=0\count142=0\count143=0\count151=0\count152=0\count153=0
\count231=5 
\count232=0 
\count233=0 
\multiply\count231 by \count211 
\multiply\count232 by \count212 
\multiply\count233 by \count213 
\advance\count141 by \count231
\advance\count141 by \count232
\advance\count141 by \count233
\divide\count141 by 5 
\count231=5 
\count232=0 
\count233=0 
\multiply\count231 by \count221 
\multiply\count232 by \count222 
\multiply\count233 by \count223 
\advance\count151 by \count231
\advance\count151 by \count232
\advance\count151 by \count233
\divide\count151 by 5
\count231=2 
\count232=2 
\count233=1 
\multiply\count231 by \count211 
\multiply\count232 by \count212 
\multiply\count233 by \count213 
\advance\count142 by \count231
\advance\count142 by \count232
\advance\count142 by \count233
\divide\count142 by 5 
\count231=2 
\count232=2 
\count233=1 
\multiply\count231 by \count221 
\multiply\count232 by \count222 
\multiply\count233 by \count223 
\advance\count152 by \count231
\advance\count152 by \count232
\advance\count152 by \count233
\divide\count152 by 5
\count231=2 
\count232=1 
\count233=2 
\multiply\count231 by \count211 
\multiply\count232 by \count212 
\multiply\count233 by \count213 
\advance\count143 by \count231
\advance\count143 by \count232
\advance\count143 by \count233
\divide\count143 by 5 
\count231=2 
\count232=1 
\count233=2 
\multiply\count231 by \count221 
\multiply\count232 by \count222 
\multiply\count233 by \count223 
\advance\count153 by \count231
\advance\count153 by \count232
\advance\count153 by \count233
\divide\count153 by 5
\or
\count141=0\count142=0\count143=0\count151=0\count152=0\count153=0
\count231=2 
\count232=2 
\count233=1 
\multiply\count231 by \count211 
\multiply\count232 by \count212 
\multiply\count233 by \count213 
\advance\count141 by \count231
\advance\count141 by \count232
\advance\count141 by \count233
\divide\count141 by 5 
\count231=2 
\count232=2 
\count233=1 
\multiply\count231 by \count221 
\multiply\count232 by \count222 
\multiply\count233 by \count223 
\advance\count151 by \count231
\advance\count151 by \count232
\advance\count151 by \count233
\divide\count151 by 5
\count231=0 
\count232=5 
\count233=0 
\multiply\count231 by \count211 
\multiply\count232 by \count212 
\multiply\count233 by \count213 
\advance\count142 by \count231
\advance\count142 by \count232
\advance\count142 by \count233
\divide\count142 by 5 
\count231=0 
\count232=5 
\count233=0 
\multiply\count231 by \count221 
\multiply\count232 by \count222 
\multiply\count233 by \count223 
\advance\count152 by \count231
\advance\count152 by \count232
\advance\count152 by \count233
\divide\count152 by 5
\count231=1 
\count232=2 
\count233=2 
\multiply\count231 by \count211 
\multiply\count232 by \count212 
\multiply\count233 by \count213 
\advance\count143 by \count231
\advance\count143 by \count232
\advance\count143 by \count233
\divide\count143 by 5 
\count231=1 
\count232=2 
\count233=2 
\multiply\count231 by \count221 
\multiply\count232 by \count222 
\multiply\count233 by \count223 
\advance\count153 by \count231
\advance\count153 by \count232
\advance\count153 by \count233
\divide\count153 by 5
\or
\count141=0\count142=0\count143=0\count151=0\count152=0\count153=0
\count231=2 
\count232=1 
\count233=2 
\multiply\count231 by \count211 
\multiply\count232 by \count212 
\multiply\count233 by \count213 
\advance\count141 by \count231
\advance\count141 by \count232
\advance\count141 by \count233
\divide\count141 by 5 
\count231=2 
\count232=1 
\count233=2 
\multiply\count231 by \count221 
\multiply\count232 by \count222 
\multiply\count233 by \count223 
\advance\count151 by \count231
\advance\count151 by \count232
\advance\count151 by \count233
\divide\count151 by 5
\count231=1 
\count232=2 
\count233=2 
\multiply\count231 by \count211 
\multiply\count232 by \count212 
\multiply\count233 by \count213 
\advance\count142 by \count231
\advance\count142 by \count232
\advance\count142 by \count233
\divide\count142 by 5 
\count231=1 
\count232=2 
\count233=2 
\multiply\count231 by \count221 
\multiply\count232 by \count222 
\multiply\count233 by \count223 
\advance\count152 by \count231
\advance\count152 by \count232
\advance\count152 by \count233
\divide\count152 by 5
\count231=0 
\count232=0 
\count233=5 
\multiply\count231 by \count211 
\multiply\count232 by \count212 
\multiply\count233 by \count213 
\advance\count143 by \count231
\advance\count143 by \count232
\advance\count143 by \count233
\divide\count143 by 5 
\count231=0 
\count232=0 
\count233=5 
\multiply\count231 by \count221 
\multiply\count232 by \count222 
\multiply\count233 by \count223 
\advance\count153 by \count231
\advance\count153 by \count232
\advance\count153 by \count233
\divide\count153 by 5
\fi%
\count211=\count141
\count212=\count142
\count213=\count143
\count221=\count151
\count222=\count152
\count223=\count153
}


\def\serr#1#2{ \power{3}{#1}
\loop 
\advance\count91 by-1 
\count92=#1 \count95=\count91 
\count211=0
\count212=62500
\count213=125000
\count221=0
\count222=109375
\count223=0
{\loop 
\advance\count92 by-1 
\count93=\count95 
\divide\count95 by 3 
\count94=\count95 
\multiply\count94 by 3 
\advance\count93 by-\count94 
\trooh{\count93} 
\ifnum \count92>0 \repeat \trigg{\count211}{\count221}{\count212}{\count222}{\count213}{\count223}{#2}}
\ifnum \count91>0 \repeat} 


\BFIG{fig-Sig}{\Sig\ in the standard harmonic coordinates.}
{\begin{picture}(195,170)(0,0)
\hhh\end{picture}}
\BEX{e-SG}{\Sig}{The \Sig\ is a \frf. The standard energy 
form \cite{Ki0,Ki1,Ki}  on the \Sig\  satisfies all the assumptions of our paper. 
The \Sig\ in \hc, see \Fig{fig-Sig}, was first considered in \cite{Ki2}, 
where the statement of 
\Thm{thm-C1} was proved in this case. The statement of
 \Thm{thm-C2} was announced in \cite{T2} without a proof.
In the case of the standard energy 
form on the \Sig\ \Conj{conj-Z} was proved in \cite{Ku1}. 
The fact that the energy measure is singular with respect to any product (Bernoulli) measure was proved in \cite{Ku1,BST,Hino1,Hino2}.}

\def\CIRCLEE#1{
\count101=#1
\count102=#1
\count103=#1
\multiply\count102 by 414
\multiply\count103 by 707
\divide\count101 by 2
\divide\count102 by 2000
\divide\count103 by 2000
\qbezier(\count101,0)(\count101,\count102)(\count103,\count103)
\qbezier(0,\count101)(\count102,\count101)(\count103,\count103)
\qbezier(-\count101,0)(-\count101,\count102)(-\count103,\count103)
\qbezier(-0,\count101)(-\count102,\count101)(-\count103,\count103)
\qbezier(-\count101,-0)(-\count101,-\count102)(-\count103,-\count103)
\qbezier(-0,-\count101)(-\count102,-\count101)(-\count103,-\count103)
\qbezier(\count101,-0)(\count101,-\count102)(\count103,-\count103)
\qbezier(0,-\count101)(\count102,-\count101)(\count103,-\count103)
}

\def\CIRCLE#1#2#3{\put(#1,#2){\CIRCLEE{#3}}
\put(-#1,#2){\CIRCLEE{#3}}
\count101=#1
\count102=#2
\count111=#1
\count112=#2
\multiply\count102 by 866
\divide\count102 by 1000
\divide\count101 by-2
\advance\count101 by \count102 
\multiply\count111 by 866
\divide\count111 by 1000
\divide\count112 by 2
\advance\count111 by \count112 
\put(\count101,-\count111){\CIRCLEE{#3}}
\put(-\count101,-\count111){\CIRCLEE{#3}}
\count101=-#1
\count102=#2
\count111=-#1
\count112=#2
\multiply\count102 by 866
\divide\count102 by 1000
\divide\count101 by-2
\advance\count101 by \count102 
\multiply\count111 by 866
\divide\count111 by 1000
\divide\count112 by 2
\advance\count111 by \count112 
\put(\count101,-\count111){\CIRCLEE{#3}}
\put(-\count101,-\count111){\CIRCLEE{#3}}
}

\BFIG{fig-Apollo}{The residue set of the Apollonian packing.}
{\begin{picture}(173.2,150)(-86.6,-50)
\hhh\end{picture}}

\BEX{e-Apollo}{The residue set of the Apollonian packing}
{It was proved in \cite{T2} that the residue set of the Apollonian packing, see \Fig{fig-Apollo},
is the \Sig\ in harmonic coordinates defined by a  non \s-s\ resistance form. This resistance form 
satisfies all the assumptions of our paper, including the (HC) assumption. }

\BEX{e-SGr}{Random \Sig s}{In \cite{MST} a family of random \Sig s was described using harmonic coordinates. Naturally, the results of this paper apply to these random gaskets, and the (HC) assumption is satisfied due to the way in which these gaskets are constructed. Also, many examples of random fractals in \cite{H,H2} satisfy the (HC) assumption, although the harmonic coordinates were not considered explicitly. }

\def\trigg#1#2#3#4#5#6#7{{\LiNe{#1}{#2}{#3}{#4}{#7} 
\LiNe{#3}{#4}{#5}{#6}{#7}\LiNe{#5}{#6}{#1}{#2}{#7}}}

\def\trig#1#2#3#4#5#6{
\count211=#3
\count212=#4
\count221=#1
\count222=#2
\count231=#1
\count232=#2
\advance\count211 by #5
\advance\count212 by #6
\advance\count221 by #5
\advance\count222 by #6
\advance\count231 by #3
\advance\count232 by #4
\count233=#1
\count234=#2
\count235=#3
\count236=#4
\multiply\count233 by 2
\multiply\count234 by 2
\multiply\count235 by 2
\multiply\count236 by 2
\qbezier(\count233,\count234)(\count231,\count232)(\count235,\count236)
\count233=#1
\count234=#2
\count235=#5
\count236=#6
\multiply\count233 by 2
\multiply\count234 by 2
\multiply\count235 by 2
\multiply\count236 by 2
\qbezier(\count233,\count234)(\count221,\count222)(\count235,\count236)
\count233=#3
\count234=#4
\count235=#5
\count236=#6
\multiply\count233 by 2
\multiply\count234 by 2
\multiply\count235 by 2
\multiply\count236 by 2
\qbezier(\count233,\count234)(\count211,\count212)(\count235,\count236)
}



\def\power#1#2{\count91=1 \count92=#2 \loop 
\advance\count92 by-1 \multiply\count91 by #1 
\ifnum \count92>0 \repeat}

\def\hrooh#1{
\ifcase#1 
\count141=0\count142=0\count143=0\count151=0\count152=0\count153=0
\count231=7 
\count232=0 
\count233=0 
\multiply\count231 by \count211 
\multiply\count232 by \count212 
\multiply\count233 by \count213 
\advance\count141 by \count231
\advance\count141 by \count232
\advance\count141 by \count233
\divide\count141 by 7 
\count231=7 
\count232=0 
\count233=0 
\multiply\count231 by \count221 
\multiply\count232 by \count222 
\multiply\count233 by \count223 
\advance\count151 by \count231
\advance\count151 by \count232
\advance\count151 by \count233
\divide\count151 by 7
\count231=4 
\count232=2 
\count233=1 
\multiply\count231 by \count211 
\multiply\count232 by \count212 
\multiply\count233 by \count213 
\advance\count142 by \count231
\advance\count142 by \count232
\advance\count142 by \count233
\divide\count142 by 7 
\count231=4 
\count232=2 
\count233=1 
\multiply\count231 by \count221 
\multiply\count232 by \count222 
\multiply\count233 by \count223 
\advance\count152 by \count231
\advance\count152 by \count232
\advance\count152 by \count233
\divide\count152 by 7
\count231=4 
\count232=1 
\count233=2 
\multiply\count231 by \count211 
\multiply\count232 by \count212 
\multiply\count233 by \count213 
\advance\count143 by \count231
\advance\count143 by \count232
\advance\count143 by \count233
\divide\count143 by 7 
\count231=4 
\count232=1 
\count233=2 
\multiply\count231 by \count221 
\multiply\count232 by \count222 
\multiply\count233 by \count223 
\advance\count153 by \count231
\advance\count153 by \count232
\advance\count153 by \count233
\divide\count153 by 7
\or
\count141=0\count142=0\count143=0\count151=0\count152=0\count153=0
\count231=2 
\count232=4 
\count233=1 
\multiply\count231 by \count211 
\multiply\count232 by \count212 
\multiply\count233 by \count213 
\advance\count141 by \count231
\advance\count141 by \count232
\advance\count141 by \count233
\divide\count141 by 7 
\count231=2 
\count232=4 
\count233=1 
\multiply\count231 by \count221 
\multiply\count232 by \count222 
\multiply\count233 by \count223 
\advance\count151 by \count231
\advance\count151 by \count232
\advance\count151 by \count233
\divide\count151 by 7
\count231=0 
\count232=7 
\count233=0 
\multiply\count231 by \count211 
\multiply\count232 by \count212 
\multiply\count233 by \count213 
\advance\count142 by \count231
\advance\count142 by \count232
\advance\count142 by \count233
\divide\count142 by 7 
\count231=0 
\count232=7 
\count233=0 
\multiply\count231 by \count221 
\multiply\count232 by \count222 
\multiply\count233 by \count223 
\advance\count152 by \count231
\advance\count152 by \count232
\advance\count152 by \count233
\divide\count152 by 7
\count231=1 
\count232=4 
\count233=2 
\multiply\count231 by \count211 
\multiply\count232 by \count212 
\multiply\count233 by \count213 
\advance\count143 by \count231
\advance\count143 by \count232
\advance\count143 by \count233
\divide\count143 by 7 
\count231=1 
\count232=4 
\count233=2 
\multiply\count231 by \count221 
\multiply\count232 by \count222 
\multiply\count233 by \count223 
\advance\count153 by \count231
\advance\count153 by \count232
\advance\count153 by \count233
\divide\count153 by 7
\or
\count141=0\count142=0\count143=0\count151=0\count152=0\count153=0
\count231=2 
\count232=1 
\count233=4 
\multiply\count231 by \count211 
\multiply\count232 by \count212 
\multiply\count233 by \count213 
\advance\count141 by \count231
\advance\count141 by \count232
\advance\count141 by \count233
\divide\count141 by 7 
\count231=2 
\count232=1 
\count233=4 
\multiply\count231 by \count221 
\multiply\count232 by \count222 
\multiply\count233 by \count223 
\advance\count151 by \count231
\advance\count151 by \count232
\advance\count151 by \count233
\divide\count151 by 7
\count231=1 
\count232=2 
\count233=4 
\multiply\count231 by \count211 
\multiply\count232 by \count212 
\multiply\count233 by \count213 
\advance\count142 by \count231
\advance\count142 by \count232
\advance\count142 by \count233
\divide\count142 by 7 
\count231=1 
\count232=2 
\count233=4 
\multiply\count231 by \count221 
\multiply\count232 by \count222 
\multiply\count233 by \count223 
\advance\count152 by \count231
\advance\count152 by \count232
\advance\count152 by \count233
\divide\count152 by 7
\count231=0 
\count232=0 
\count233=7 
\multiply\count231 by \count211 
\multiply\count232 by \count212 
\multiply\count233 by \count213 
\advance\count143 by \count231
\advance\count143 by \count232
\advance\count143 by \count233
\divide\count143 by 7 
\count231=0 
\count232=0 
\count233=7 
\multiply\count231 by \count221 
\multiply\count232 by \count222 
\multiply\count233 by \count223 
\advance\count153 by \count231
\advance\count153 by \count232
\advance\count153 by \count233
\divide\count153 by 7
\or
%
%
%
\count141=0\count142=0\count143=0\count151=0\count152=0\count153=0
\count231=1 
\count232=3 
\count233=3 
\multiply\count231 by \count211 
\multiply\count232 by \count212 
\multiply\count233 by \count213 
\advance\count141 by \count231
\advance\count141 by \count232
\advance\count141 by \count233
\divide\count141 by 7 
\count231=1 
\count232=3 
\count233=3 
\multiply\count231 by \count221 
\multiply\count232 by \count222 
\multiply\count233 by \count223 
\advance\count151 by \count231
\advance\count151 by \count232
\advance\count151 by \count233
\divide\count151 by 7
\count231=1 
\count232=2 
\count233=4 
\multiply\count231 by \count211 
\multiply\count232 by \count212 
\multiply\count233 by \count213 
\advance\count142 by \count231
\advance\count142 by \count232
\advance\count142 by \count233
\divide\count142 by 7 
\count231=1 
\count232=2 
\count233=4 
\multiply\count231 by \count221 
\multiply\count232 by \count222 
\multiply\count233 by \count223 
\advance\count152 by \count231
\advance\count152 by \count232
\advance\count152 by \count233
\divide\count152 by 7
\count231=1 
\count232=4 
\count233=2 
\multiply\count231 by \count211 
\multiply\count232 by \count212 
\multiply\count233 by \count213 
\advance\count143 by \count231
\advance\count143 by \count232
\advance\count143 by \count233
\divide\count143 by 7 
\count231=1 
\count232=4 
\count233=2 
\multiply\count231 by \count221 
\multiply\count232 by \count222 
\multiply\count233 by \count223 
\advance\count153 by \count231
\advance\count153 by \count232
\advance\count153 by \count233
\divide\count153 by 7
\or
\count141=0\count142=0\count143=0\count151=0\count152=0\count153=0
\count231=2 
\count232=1 
\count233=4 
\multiply\count231 by \count211 
\multiply\count232 by \count212 
\multiply\count233 by \count213 
\advance\count141 by \count231
\advance\count141 by \count232
\advance\count141 by \count233
\divide\count141 by 7 
\count231=2 
\count232=1 
\count233=4 
\multiply\count231 by \count221 
\multiply\count232 by \count222 
\multiply\count233 by \count223 
\advance\count151 by \count231
\advance\count151 by \count232
\advance\count151 by \count233
\divide\count151 by 7
\count231=3 
\count232=1 
\count233=3 
\multiply\count231 by \count211 
\multiply\count232 by \count212 
\multiply\count233 by \count213 
\advance\count142 by \count231
\advance\count142 by \count232
\advance\count142 by \count233
\divide\count142 by 7 
\count231=3 
\count232=1 
\count233=3 
\multiply\count231 by \count221 
\multiply\count232 by \count222 
\multiply\count233 by \count223 
\advance\count152 by \count231
\advance\count152 by \count232
\advance\count152 by \count233
\divide\count152 by 7
\count231=4 
\count232=1 
\count233=2 
\multiply\count231 by \count211 
\multiply\count232 by \count212 
\multiply\count233 by \count213 
\advance\count143 by \count231
\advance\count143 by \count232
\advance\count143 by \count233
\divide\count143 by 7 
\count231=4 
\count232=1 
\count233=2 
\multiply\count231 by \count221 
\multiply\count232 by \count222 
\multiply\count233 by \count223 
\advance\count153 by \count231
\advance\count153 by \count232
\advance\count153 by \count233
\divide\count153 by 7
\or
\count141=0\count142=0\count143=0\count151=0\count152=0\count153=0
\count231=2 
\count232=4 
\count233=1 
\multiply\count231 by \count211 
\multiply\count232 by \count212 
\multiply\count233 by \count213 
\advance\count141 by \count231
\advance\count141 by \count232
\advance\count141 by \count233
\divide\count141 by 7 
\count231=2 
\count232=4 
\count233=1 
\multiply\count231 by \count221 
\multiply\count232 by \count222 
\multiply\count233 by \count223 
\advance\count151 by \count231
\advance\count151 by \count232
\advance\count151 by \count233
\divide\count151 by 7
\count231=4 
\count232=2 
\count233=1 
\multiply\count231 by \count211 
\multiply\count232 by \count212 
\multiply\count233 by \count213 
\advance\count142 by \count231
\advance\count142 by \count232
\advance\count142 by \count233
\divide\count142 by 7 
\count231=4 
\count232=2 
\count233=1 
\multiply\count231 by \count221 
\multiply\count232 by \count222 
\multiply\count233 by \count223 
\advance\count152 by \count231
\advance\count152 by \count232
\advance\count152 by \count233
\divide\count152 by 7
\count231=3 
\count232=3 
\count233=1 
\multiply\count231 by \count211 
\multiply\count232 by \count212 
\multiply\count233 by \count213 
\advance\count143 by \count231
\advance\count143 by \count232
\advance\count143 by \count233
\divide\count143 by 7 
\count231=3 
\count232=3 
\count233=1 
\multiply\count231 by \count221 
\multiply\count232 by \count222 
\multiply\count233 by \count223 
\advance\count153 by \count231
\advance\count153 by \count232
\advance\count153 by \count233
\divide\count153 by 7
\fi%
\count211=\count141
\count212=\count142
\count213=\count143
\count221=\count151
\count222=\count152
\count223=\count153
}


\def\herr#1#2{ \power{6}{#1}
\loop 
\advance\count91 by-1 
\count92=#1 \count95=\count91 
\count211=0
\count212=62500
\count213=125000
\count221=0
\count222=109375
\count223=0
{\loop 
\advance\count92 by-1 
\count93=\count95 
\divide\count95 by 6 
\count94=\count95 
\multiply\count94 by 6 
\advance\count93 by-\count94 
\hrooh{\count93} 
\ifnum \count92>0 
\repeat \trigg{\count211}{\count221}{\count212}{\count222}{\count213}{\count223}{#2}}
\ifnum \count91>0 
\repeat} 

\BFIG{fig-HexH}{The  hexagasket in harmonic coordinates and its  first approximation.}
{\begin{picture}(195,170)(0,0) 
\hhh\end{picture}

\smallskip

\begin{picture}(210,220)(-105,-5)
\hhh\end{picture}}

\BEX{e-X}{Hexagasket}{According to   \cite{T1}, the Hexagasket satisfies the (WN) assumption but not the (HC) assumption. However, by small perturbations of the harmonic coordinates  one can  construct two functions of finite energy which 
 map the hexagasket into $\mathbb R^2$ homeomorphically. Then the conclusion of Theorems \ref{thm-Z} and \ref{thm-C1} will hold because of the general theory of Dirichlet forms in \cite{BouH,FOT} (see \Rem{rem-Z-dnu}). However \Thm{thm-C2} will not hold unless these coordinates are in the domain of the energy \Lp, which is difficult to verify. }


\BEX{e-Qpcf}{Quotients of \pcf\ fractals}{If we 
consider quotient of a \pcf\ fractal defined by its space of harmonic functions, and conditions of \Thm{cor-HC} are satisfied (see also \Thm{prop-supp-nu}), then we have a \frf\ which 
satisfies the (HC) assumption by definition. In the case of the 
Hexagasket this is illustrated in Figure~\ref{fig-HexH}. Note that this set is not self-affine. In harmonic coordinates the Hexagasket is represented as a union of a Cantor set and a disjoint union of countably many closed straight line intervals. 
One can show that the energy measure of this Cantor set is zero, and in fact the energy measure is proportional to the Lebesgue measure on each segment. 
Note that in the limit no two intervals meet and so it is not a quantum graph, but can be called a generalized quantum graph. In this case a three point boundary, see \cite{St4}, is chosen so that the resulting fractal can be embedded in $\mathbb R^2$. For a different choice of the boundary the local structure of the fractal in \hc\ is the same. }

\BEX{e-V}{Vicsek set}{Vicsek set (see, for instance, \cite{T2}) is a \frf\ which does not satisfy 
 the (WN) and (HC) assumptions. 
In harmonic coordinates it is represented by four straight line segments which are joined at a point. Therefore in our construction $F_H$ is a 
quantum graph with five vertices and four edges, which is not homeomorphic to the Vicsek set.}

\def\pciSigoo#1{
\ifcase#1 
%
%
%
\count141=0\count142=0\count143=0\count151=0\count152=0\count153=0
\count231=159 
\count232=0 
\count233=0 
\multiply\count231 by \count211 
\multiply\count232 by \count212 
\multiply\count233 by \count213 
\advance\count141 by \count231
\advance\count141 by \count232
\advance\count141 by \count233
\divide\count141 by 159 
\count231=159 
\count232=0 
\count233=0 
\multiply\count231 by \count221 
\multiply\count232 by \count222 
\multiply\count233 by \count223 
\advance\count151 by \count231
\advance\count151 by \count232
\advance\count151 by \count233
\divide\count151 by 159
\count231=95 
\count232=38 
\count233=26 
\multiply\count231 by \count211 
\multiply\count232 by \count212 
\multiply\count233 by \count213 
\advance\count142 by \count231
\advance\count142 by \count232
\advance\count142 by \count233
\divide\count142 by 159 
\count231=95 
\count232=38 
\count233=26 
\multiply\count231 by \count221 
\multiply\count232 by \count222 
\multiply\count233 by \count223 
\advance\count152 by \count231
\advance\count152 by \count232
\advance\count152 by \count233
\divide\count152 by 159
\count231=73 
\count232=73 
\count233=13 
\multiply\count231 by \count211 
\multiply\count232 by \count212 
\multiply\count233 by \count213 
\advance\count143 by \count231
\advance\count143 by \count232
\advance\count143 by \count233
\divide\count143 by 159 
\count231=73 
\count232=73 
\count233=13 
\multiply\count231 by \count221 
\multiply\count232 by \count222 
\multiply\count233 by \count223 
\advance\count153 by \count231
\advance\count153 by \count232
\advance\count153 by \count233
\divide\count153 by 159
\or
\count141=0\count142=0\count143=0\count151=0\count152=0\count153=0
\count231=159 
\count232=0 
\count233=0 
\multiply\count231 by \count211 
\multiply\count232 by \count212 
\multiply\count233 by \count213 
\advance\count141 by \count231
\advance\count141 by \count232
\advance\count141 by \count233
\divide\count141 by 159 
\count231=159 
\count232=0 
\count233=0 
\multiply\count231 by \count221 
\multiply\count232 by \count222 
\multiply\count233 by \count223 
\advance\count151 by \count231
\advance\count151 by \count232
\advance\count151 by \count233
\divide\count151 by 159
\count231=95 
\count232=26 
\count233=38 
\multiply\count231 by \count211 
\multiply\count232 by \count212 
\multiply\count233 by \count213 
\advance\count142 by \count231
\advance\count142 by \count232
\advance\count142 by \count233
\divide\count142 by 159 
\count231=95 
\count232=26 
\count233=38 
\multiply\count231 by \count221 
\multiply\count232 by \count222 
\multiply\count233 by \count223 
\advance\count152 by \count231
\advance\count152 by \count232
\advance\count152 by \count233
\divide\count152 by 159
\count231=73 
\count232=13 
\count233=73 
\multiply\count231 by \count211 
\multiply\count232 by \count212 
\multiply\count233 by \count213 
\advance\count143 by \count231
\advance\count143 by \count232
\advance\count143 by \count233
\divide\count143 by 159 
\count231=73 
\count232=13 
\count233=73 
\multiply\count231 by \count221 
\multiply\count232 by \count222 
\multiply\count233 by \count223 
\advance\count153 by \count231
\advance\count153 by \count232
\advance\count153 by \count233
\divide\count153 by 159
\or
\count141=0\count142=0\count143=0\count151=0\count152=0\count153=0
\count231=13 
\count232=73 
\count233=73 
\multiply\count231 by \count211 
\multiply\count232 by \count212 
\multiply\count233 by \count213 
\advance\count141 by \count231
\advance\count141 by \count232
\advance\count141 by \count233
\divide\count141 by 159 
\count231=13 
\count232=73 
\count233=73 
\multiply\count231 by \count221 
\multiply\count232 by \count222 
\multiply\count233 by \count223 
\advance\count151 by \count231
\advance\count151 by \count232
\advance\count151 by \count233
\divide\count151 by 159
\count231=0
\count232=159
\count233=0
\multiply\count231 by \count211 
\multiply\count232 by \count212 
\multiply\count233 by \count213 
\advance\count142 by \count231
\advance\count142 by \count232
\advance\count142 by \count233
\divide\count142 by 159 
\count231=0
\count232=159
\count233=0
\multiply\count231 by \count221 
\multiply\count232 by \count222 
\multiply\count233 by \count223 
\advance\count152 by \count231
\advance\count152 by \count232
\advance\count152 by \count233
\divide\count152 by 159
\count231=26 
\count232=95 
\count233=38 
\multiply\count231 by \count211 
\multiply\count232 by \count212 
\multiply\count233 by \count213 
\advance\count143 by \count231
\advance\count143 by \count232
\advance\count143 by \count233
\divide\count143 by 159 
\count231=26 
\count232=95 
\count233=38 
\multiply\count231 by \count221 
\multiply\count232 by \count222 
\multiply\count233 by \count223 
\advance\count153 by \count231
\advance\count153 by \count232
\advance\count153 by \count233
\divide\count153 by 159
\or
\count141=0\count142=0\count143=0\count151=0\count152=0\count153=0
\count231=73 
\count232=73 
\count233=13 
\multiply\count231 by \count211 
\multiply\count232 by \count212 
\multiply\count233 by \count213 
\advance\count141 by \count231
\advance\count141 by \count232
\advance\count141 by \count233
\divide\count141 by 159 
\count231=73 
\count232=73 
\count233=13 
\multiply\count231 by \count221 
\multiply\count232 by \count222 
\multiply\count233 by \count223 
\advance\count151 by \count231
\advance\count151 by \count232
\advance\count151 by \count233
\divide\count151 by 159
\count231=0
\count232=159
\count233=0
\multiply\count231 by \count211 
\multiply\count232 by \count212 
\multiply\count233 by \count213 
\advance\count142 by \count231
\advance\count142 by \count232
\advance\count142 by \count233
\divide\count142 by 159 
\count231=0
\count232=159
\count233=0
\multiply\count231 by \count221 
\multiply\count232 by \count222 
\multiply\count233 by \count223 
\advance\count152 by \count231
\advance\count152 by \count232
\advance\count152 by \count233
\divide\count152 by 159
\count231=38 
\count232=95 
\count233=26 
\multiply\count231 by \count211 
\multiply\count232 by \count212 
\multiply\count233 by \count213 
\advance\count143 by \count231
\advance\count143 by \count232
\advance\count143 by \count233
\divide\count143 by 159 
\count231=38 
\count232=95 
\count233=26 
\multiply\count231 by \count221 
\multiply\count232 by \count222 
\multiply\count233 by \count223 
\advance\count153 by \count231
\advance\count153 by \count232
\advance\count153 by \count233
\divide\count153 by 159
\or
\count141=0\count142=0\count143=0\count151=0\count152=0\count153=0
\count231=26
\count232=38
\count233=95
\multiply\count231 by \count211 
\multiply\count232 by \count212 
\multiply\count233 by \count213 
\advance\count141 by \count231
\advance\count141 by \count232
\advance\count141 by \count233
\divide\count141 by 159 
\count231=26 
\count232=38 
\count233=95 
\multiply\count231 by \count221 
\multiply\count232 by \count222 
\multiply\count233 by \count223 
\advance\count151 by \count231
\advance\count151 by \count232
\advance\count151 by \count233
\divide\count151 by 159
\count231=13
\count232=73
\count233=73
\multiply\count231 by \count211 
\multiply\count232 by \count212 
\multiply\count233 by \count213 
\advance\count142 by \count231
\advance\count142 by \count232
\advance\count142 by \count233
\divide\count142 by 159 
\count231=13
\count232=73
\count233=73
\multiply\count231 by \count221 
\multiply\count232 by \count222 
\multiply\count233 by \count223 
\advance\count152 by \count231
\advance\count152 by \count232
\advance\count152 by \count233
\divide\count152 by 159
\count231=0 
\count232=0 
\count233=159 
\multiply\count231 by \count211 
\multiply\count232 by \count212 
\multiply\count233 by \count213 
\advance\count143 by \count231
\advance\count143 by \count232
\advance\count143 by \count233
\divide\count143 by 159 
\count231=0 
\count232=0 
\count233=159 
\multiply\count231 by \count221 
\multiply\count232 by \count222 
\multiply\count233 by \count223 
\advance\count153 by \count231
\advance\count153 by \count232
\advance\count153 by \count233
\divide\count153 by 159
\or
\count141=0\count142=0\count143=0\count151=0\count152=0\count153=0
\count231=38
\count232=26
\count233=95
\multiply\count231 by \count211 
\multiply\count232 by \count212 
\multiply\count233 by \count213 
\advance\count141 by \count231
\advance\count141 by \count232
\advance\count141 by \count233
\divide\count141 by 159 
\count231=38 
\count232=26 
\count233=95 
\multiply\count231 by \count221 
\multiply\count232 by \count222 
\multiply\count233 by \count223 
\advance\count151 by \count231
\advance\count151 by \count232
\advance\count151 by \count233
\divide\count151 by 159
\count231=73
\count232=13
\count233=73
\multiply\count231 by \count211 
\multiply\count232 by \count212 
\multiply\count233 by \count213 
\advance\count142 by \count231
\advance\count142 by \count232
\advance\count142 by \count233
\divide\count142 by 159 
\count231=73
\count232=13
\count233=73
\multiply\count231 by \count221 
\multiply\count232 by \count222 
\multiply\count233 by \count223 
\advance\count152 by \count231
\advance\count152 by \count232
\advance\count152 by \count233
\divide\count152 by 159
\count231=0 
\count232=0 
\count233=159 
\multiply\count231 by \count211 
\multiply\count232 by \count212 
\multiply\count233 by \count213 
\advance\count143 by \count231
\advance\count143 by \count232
\advance\count143 by \count233
\divide\count143 by 159 
\count231=0 
\count232=0 
\count233=159 
\multiply\count231 by \count221 
\multiply\count232 by \count222 
\multiply\count233 by \count223 
\advance\count153 by \count231
\advance\count153 by \count232
\advance\count153 by \count233
\divide\count153 by 159
\or
\count141=0\count142=0\count143=0\count151=0\count152=0\count153=0
\count231=38
\count232=26
\count233=95
\multiply\count231 by \count211 
\multiply\count232 by \count212 
\multiply\count233 by \count213 
\advance\count141 by \count231
\advance\count141 by \count232
\advance\count141 by \count233
\divide\count141 by 159 
\count231=38 
\count232=26 
\count233=95 
\multiply\count231 by \count221 
\multiply\count232 by \count222 
\multiply\count233 by \count223 
\advance\count151 by \count231
\advance\count151 by \count232
\advance\count151 by \count233
\divide\count151 by 159
\count231=26
\count232=38
\count233=95
\multiply\count231 by \count211 
\multiply\count232 by \count212 
\multiply\count233 by \count213 
\advance\count142 by \count231
\advance\count142 by \count232
\advance\count142 by \count233
\divide\count142 by 159 
\count231=26
\count232=38
\count233=95
\multiply\count231 by \count221 
\multiply\count232 by \count222 
\multiply\count233 by \count223 
\advance\count152 by \count231
\advance\count152 by \count232
\advance\count152 by \count233
\divide\count152 by 159
\count231=53 
\count232=53 
\count233=53 
\multiply\count231 by \count211 
\multiply\count232 by \count212 
\multiply\count233 by \count213 
\advance\count143 by \count231
\advance\count143 by \count232
\advance\count143 by \count233
\divide\count143 by 159 
\count231=53 
\count232=53 
\count233=53 
\multiply\count231 by \count221 
\multiply\count232 by \count222 
\multiply\count233 by \count223 
\advance\count153 by \count231
\advance\count153 by \count232
\advance\count153 by \count233
\divide\count153 by 159
\or
\count141=0\count142=0\count143=0\count151=0\count152=0\count153=0
\count231=26
\count232=95
\count233=38
\multiply\count231 by \count211 
\multiply\count232 by \count212 
\multiply\count233 by \count213 
\advance\count141 by \count231
\advance\count141 by \count232
\advance\count141 by \count233
\divide\count141 by 159 
\count231=26 
\count232=95 
\count233=38 
\multiply\count231 by \count221 
\multiply\count232 by \count222 
\multiply\count233 by \count223 
\advance\count151 by \count231
\advance\count151 by \count232
\advance\count151 by \count233
\divide\count151 by 159
\count231=53
\count232=53
\count233=53
\multiply\count231 by \count211 
\multiply\count232 by \count212 
\multiply\count233 by \count213 
\advance\count142 by \count231
\advance\count142 by \count232
\advance\count142 by \count233
\divide\count142 by 159 
\count231=53
\count232=53
\count233=53
\multiply\count231 by \count221 
\multiply\count232 by \count222 
\multiply\count233 by \count223 
\advance\count152 by \count231
\advance\count152 by \count232
\advance\count152 by \count233
\divide\count152 by 159
\count231=38 
\count232=95 
\count233=26 
\multiply\count231 by \count211 
\multiply\count232 by \count212 
\multiply\count233 by \count213 
\advance\count143 by \count231
\advance\count143 by \count232
\advance\count143 by \count233
\divide\count143 by 159 
\count231=38 
\count232=95 
\count233=26 
\multiply\count231 by \count221 
\multiply\count232 by \count222 
\multiply\count233 by \count223 
\advance\count153 by \count231
\advance\count153 by \count232
\advance\count153 by \count233
\divide\count153 by 159
\or
\count141=0\count142=0\count143=0\count151=0\count152=0\count153=0
\count231=53
\count232=53
\count233=53
\multiply\count231 by \count211 
\multiply\count232 by \count212 
\multiply\count233 by \count213 
\advance\count141 by \count231
\advance\count141 by \count232
\advance\count141 by \count233
\divide\count141 by 159 
\count231=53 
\count232=53 
\count233=53 
\multiply\count231 by \count221 
\multiply\count232 by \count222 
\multiply\count233 by \count223 
\advance\count151 by \count231
\advance\count151 by \count232
\advance\count151 by \count233
\divide\count151 by 159
\count231=95
\count232=38
\count233=26
\multiply\count231 by \count211 
\multiply\count232 by \count212 
\multiply\count233 by \count213 
\advance\count142 by \count231
\advance\count142 by \count232
\advance\count142 by \count233
\divide\count142 by 159 
\count231=95
\count232=38
\count233=26
\multiply\count231 by \count221 
\multiply\count232 by \count222 
\multiply\count233 by \count223 
\advance\count152 by \count231
\advance\count152 by \count232
\advance\count152 by \count233
\divide\count152 by 159
\count231=95 
\count232=26 
\count233=38 
\multiply\count231 by \count211 
\multiply\count232 by \count212 
\multiply\count233 by \count213 
\advance\count143 by \count231
\advance\count143 by \count232
\advance\count143 by \count233
\divide\count143 by 159 
\count231=95 
\count232=26 
\count233=38 
\multiply\count231 by \count221 
\multiply\count232 by \count222 
\multiply\count233 by \count223 
\advance\count153 by \count231
\advance\count153 by \count232
\advance\count153 by \count233
\divide\count153 by 159
\fi%
\count211=\count141
\count212=\count142
\count213=\count143
\count221=\count151
\count222=\count152
\count223=\count153
}

\def\trigg#1#2#3#4#5#6#7{{
\count107=#1 
\count108=#2 
\advance\count107 by #3 \advance\count107 by #5 \divide\count107 by 3 
\advance\count108 by #4 \advance\count108 by #6 \divide\count108 by 3
\LiNe{#1}{#2}{\count107}{\count108}{#7} 
\LiNe{#3}{#4}{\count107}{\count108}{#7} 
\LiNe{#5}{#6}{\count107}{\count108}{#7} 
}}

\def\trigg#1#2#3#4#5#6#7{{\LiNe{#1}{#2}{#3}{#4}{#7} 
\LiNe{#3}{#4}{#5}{#6}{#7}\LiNe{#5}{#6}{#1}{#2}{#7}}}

\def\power#1#2{\count91=1 \count92=#2 \loop 
\advance\count92 by-1 \multiply\count91 by #1 
\ifnum \count92>0 \repeat}


\def\pciSig#1#2{ \power{9}{#1}
\loop 
\advance\count91 by-1 
\count92=#1 \count95=\count91 
\count211=0
\count212=62500
\count213=125000
\count221=0
\count222=109375
\count223=0
{\loop 
\advance\count92 by-1 
\count93=\count95 
\divide\count95 by 9 
\count94=\count95 
\multiply\count94 by 9 
\advance\count93 by-\count94 
\pciSigoo{\count93} 
\ifnum \count92>0 
\repeat \trigg{\count211}{\count221}{\count212}{\count222}{\count213}{\count223}{#2}}
\ifnum \count91>0 
\repeat} 

%
\def\trig#1#2#3#4#5#6{
\count211=#3
\count212=#4
\count221=#1
\count222=#2
\count231=#1
\count232=#2
\advance\count211 by #5
\divide\count211 by 2
\advance\count212 by #6
\divide\count212 by 2
\advance\count221 by #5
\divide\count221 by 2
\advance\count222 by #6
\divide\count222 by 2
\advance\count231 by #3
\divide\count231 by 2
\advance\count232 by #4
\divide\count232 by 2
\qbezier(#1,#2)(\count231,\count232)(#3,#4)
\qbezier(#1,#2)(\count221,\count222)(#5,#6)
\qbezier(#3,#4)(\count211,\count212)(#5,#6)
}
\def\trig#1#2#3#4#5#6{
\count211=#3
\count212=#4
\count221=#1
\count222=#2
\count231=#1
\count232=#2
\advance\count211 by #5
\advance\count212 by #6
\advance\count221 by #5
\advance\count222 by #6
\advance\count231 by #3
\advance\count232 by #4
\count233=#1
\count234=#2
\count235=#3
\count236=#4
\multiply\count233 by 2
\multiply\count234 by 2
\multiply\count235 by 2
\multiply\count236 by 2
\qbezier(\count233,\count234)(\count231,\count232)(\count235,\count236)
\count233=#1
\count234=#2
\count235=#5
\count236=#6
\multiply\count233 by 2
\multiply\count234 by 2
\multiply\count235 by 2
\multiply\count236 by 2
\qbezier(\count233,\count234)(\count221,\count222)(\count235,\count236)
\count233=#3
\count234=#4
\count235=#5
\count236=#6
\multiply\count233 by 2
\multiply\count234 by 2
\multiply\count235 by 2
\multiply\count236 by 2
\qbezier(\count233,\count234)(\count211,\count212)(\count235,\count236)
}

\BFIG{figpciSig}{The post-critically infinite \Sig\ in harmonic coordinates and its first approximation.}
{\begin{picture}(250,219)(0,0)
\hhh\end{picture}

\bigskip

\begin{picture}(318,245)(-159,-5)
\hhh\end{picture}}


\def\pciSigooHanother#1{
\ifcase#1 
%
%
%
\count141=0\count142=0\count143=0\count151=0\count152=0\count153=0
\count231=15 
\count232=0 
\count233=0 
\multiply\count231 by \count211 
\multiply\count232 by \count212 
\multiply\count233 by \count213 
\advance\count141 by \count231
\advance\count141 by \count232
\advance\count141 by \count233
\divide\count141 by 15 
\count231=15 
\count232=0 
\count233=0 
\multiply\count231 by \count221 
\multiply\count232 by \count222 
\multiply\count233 by \count223 
\advance\count151 by \count231
\advance\count151 by \count232
\advance\count151 by \count233
\divide\count151 by 15
\count231=11 
\count232=2 
\count233=2 
\multiply\count231 by \count211 
\multiply\count232 by \count212 
\multiply\count233 by \count213 
\advance\count142 by \count231
\advance\count142 by \count232
\advance\count142 by \count233
\divide\count142 by 15 
\count231=11 
\count232=2 
\count233=2 
\multiply\count231 by \count221 
\multiply\count232 by \count222 
\multiply\count233 by \count223 
\advance\count152 by \count231
\advance\count152 by \count232
\advance\count152 by \count233
\divide\count152 by 15
\count231=7 
\count232=7 
\count233=1 
\multiply\count231 by \count211 
\multiply\count232 by \count212 
\multiply\count233 by \count213 
\advance\count143 by \count231
\advance\count143 by \count232
\advance\count143 by \count233
\divide\count143 by 15 
\count231=7 
\count232=7 
\count233=1 
\multiply\count231 by \count221 
\multiply\count232 by \count222 
\multiply\count233 by \count223 
\advance\count153 by \count231
\advance\count153 by \count232
\advance\count153 by \count233
\divide\count153 by 15
\or
\count141=0\count142=0\count143=0\count151=0\count152=0\count153=0
\count231=15 
\count232=0 
\count233=0 
\multiply\count231 by \count211 
\multiply\count232 by \count212 
\multiply\count233 by \count213 
\advance\count141 by \count231
\advance\count141 by \count232
\advance\count141 by \count233
\divide\count141 by 15 
\count231=15 
\count232=0 
\count233=0 
\multiply\count231 by \count221 
\multiply\count232 by \count222 
\multiply\count233 by \count223 
\advance\count151 by \count231
\advance\count151 by \count232
\advance\count151 by \count233
\divide\count151 by 15
\count231=11 
\count232=2 
\count233=2 
\multiply\count231 by \count211 
\multiply\count232 by \count212 
\multiply\count233 by \count213 
\advance\count142 by \count231
\advance\count142 by \count232
\advance\count142 by \count233
\divide\count142 by 15 
\count231=11 
\count232=2 
\count233=2 
\multiply\count231 by \count221 
\multiply\count232 by \count222 
\multiply\count233 by \count223 
\advance\count152 by \count231
\advance\count152 by \count232
\advance\count152 by \count233
\divide\count152 by 15
\count231=7 
\count232=1 
\count233=7 
\multiply\count231 by \count211 
\multiply\count232 by \count212 
\multiply\count233 by \count213 
\advance\count143 by \count231
\advance\count143 by \count232
\advance\count143 by \count233
\divide\count143 by 15 
\count231=7 
\count232=1 
\count233=7 
\multiply\count231 by \count221 
\multiply\count232 by \count222 
\multiply\count233 by \count223 
\advance\count153 by \count231
\advance\count153 by \count232
\advance\count153 by \count233
\divide\count153 by 15
\or
\count141=0\count142=0\count143=0\count151=0\count152=0\count153=0
\count231=1 
\count232=7 
\count233=7 
\multiply\count231 by \count211 
\multiply\count232 by \count212 
\multiply\count233 by \count213 
\advance\count141 by \count231
\advance\count141 by \count232
\advance\count141 by \count233
\divide\count141 by 15 
\count231=1 
\count232=7 
\count233=7 
\multiply\count231 by \count221 
\multiply\count232 by \count222 
\multiply\count233 by \count223 
\advance\count151 by \count231
\advance\count151 by \count232
\advance\count151 by \count233
\divide\count151 by 15
\count231=0
\count232=15
\count233=0
\multiply\count231 by \count211 
\multiply\count232 by \count212 
\multiply\count233 by \count213 
\advance\count142 by \count231
\advance\count142 by \count232
\advance\count142 by \count233
\divide\count142 by 15 
\count231=0
\count232=15
\count233=0
\multiply\count231 by \count221 
\multiply\count232 by \count222 
\multiply\count233 by \count223 
\advance\count152 by \count231
\advance\count152 by \count232
\advance\count152 by \count233
\divide\count152 by 15
\count231=2 
\count232=11 
\count233=2 
\multiply\count231 by \count211 
\multiply\count232 by \count212 
\multiply\count233 by \count213 
\advance\count143 by \count231
\advance\count143 by \count232
\advance\count143 by \count233
\divide\count143 by 15 
\count231=2 
\count232=11 
\count233=2 
\multiply\count231 by \count221 
\multiply\count232 by \count222 
\multiply\count233 by \count223 
\advance\count153 by \count231
\advance\count153 by \count232
\advance\count153 by \count233
\divide\count153 by 15
\or
\count141=0\count142=0\count143=0\count151=0\count152=0\count153=0
\count231=7 
\count232=7 
\count233=1 
\multiply\count231 by \count211 
\multiply\count232 by \count212 
\multiply\count233 by \count213 
\advance\count141 by \count231
\advance\count141 by \count232
\advance\count141 by \count233
\divide\count141 by 15 
\count231=7 
\count232=7 
\count233=1 
\multiply\count231 by \count221 
\multiply\count232 by \count222 
\multiply\count233 by \count223 
\advance\count151 by \count231
\advance\count151 by \count232
\advance\count151 by \count233
\divide\count151 by 15
\count231=0
\count232=15
\count233=0
\multiply\count231 by \count211 
\multiply\count232 by \count212 
\multiply\count233 by \count213 
\advance\count142 by \count231
\advance\count142 by \count232
\advance\count142 by \count233
\divide\count142 by 15 
\count231=0
\count232=15
\count233=0
\multiply\count231 by \count221 
\multiply\count232 by \count222 
\multiply\count233 by \count223 
\advance\count152 by \count231
\advance\count152 by \count232
\advance\count152 by \count233
\divide\count152 by 15
\count231=2 
\count232=11 
\count233=2 
\multiply\count231 by \count211 
\multiply\count232 by \count212 
\multiply\count233 by \count213 
\advance\count143 by \count231
\advance\count143 by \count232
\advance\count143 by \count233
\divide\count143 by 15 
\count231=2 
\count232=11 
\count233=2 
\multiply\count231 by \count221 
\multiply\count232 by \count222 
\multiply\count233 by \count223 
\advance\count153 by \count231
\advance\count153 by \count232
\advance\count153 by \count233
\divide\count153 by 15
\or
\count141=0\count142=0\count143=0\count151=0\count152=0\count153=0
\count231=2
\count232=2
\count233=11
\multiply\count231 by \count211 
\multiply\count232 by \count212 
\multiply\count233 by \count213 
\advance\count141 by \count231
\advance\count141 by \count232
\advance\count141 by \count233
\divide\count141 by 15 
\count231=2 
\count232=2 
\count233=11 
\multiply\count231 by \count221 
\multiply\count232 by \count222 
\multiply\count233 by \count223 
\advance\count151 by \count231
\advance\count151 by \count232
\advance\count151 by \count233
\divide\count151 by 15
\count231=1
\count232=7
\count233=7
\multiply\count231 by \count211 
\multiply\count232 by \count212 
\multiply\count233 by \count213 
\advance\count142 by \count231
\advance\count142 by \count232
\advance\count142 by \count233
\divide\count142 by 15 
\count231=1
\count232=7
\count233=7
\multiply\count231 by \count221 
\multiply\count232 by \count222 
\multiply\count233 by \count223 
\advance\count152 by \count231
\advance\count152 by \count232
\advance\count152 by \count233
\divide\count152 by 15
\count231=0 
\count232=0 
\count233=15 
\multiply\count231 by \count211 
\multiply\count232 by \count212 
\multiply\count233 by \count213 
\advance\count143 by \count231
\advance\count143 by \count232
\advance\count143 by \count233
\divide\count143 by 15 
\count231=0 
\count232=0 
\count233=15 
\multiply\count231 by \count221 
\multiply\count232 by \count222 
\multiply\count233 by \count223 
\advance\count153 by \count231
\advance\count153 by \count232
\advance\count153 by \count233
\divide\count153 by 15
\or
\count141=0\count142=0\count143=0\count151=0\count152=0\count153=0
\count231=2
\count232=2
\count233=11
\multiply\count231 by \count211 
\multiply\count232 by \count212 
\multiply\count233 by \count213 
\advance\count141 by \count231
\advance\count141 by \count232
\advance\count141 by \count233
\divide\count141 by 15 
\count231=2 
\count232=2 
\count233=11 
\multiply\count231 by \count221 
\multiply\count232 by \count222 
\multiply\count233 by \count223 
\advance\count151 by \count231
\advance\count151 by \count232
\advance\count151 by \count233
\divide\count151 by 15
\count231=7
\count232=1
\count233=7
\multiply\count231 by \count211 
\multiply\count232 by \count212 
\multiply\count233 by \count213 
\advance\count142 by \count231
\advance\count142 by \count232
\advance\count142 by \count233
\divide\count142 by 15 
\count231=7
\count232=1
\count233=7
\multiply\count231 by \count221 
\multiply\count232 by \count222 
\multiply\count233 by \count223 
\advance\count152 by \count231
\advance\count152 by \count232
\advance\count152 by \count233
\divide\count152 by 15
\count231=0 
\count232=0 
\count233=15 
\multiply\count231 by \count211 
\multiply\count232 by \count212 
\multiply\count233 by \count213 
\advance\count143 by \count231
\advance\count143 by \count232
\advance\count143 by \count233
\divide\count143 by 15 
\count231=0 
\count232=0 
\count233=15 
\multiply\count231 by \count221 
\multiply\count232 by \count222 
\multiply\count233 by \count223 
\advance\count153 by \count231
\advance\count153 by \count232
\advance\count153 by \count233
\divide\count153 by 15
\fi%
\count211=\count141
\count212=\count142
\count213=\count143
\count221=\count151
\count222=\count152
\count223=\count153
}

\def\trigg#1#2#3#4#5#6#7{{
\count107=#1 
\count108=#2 
\advance\count107 by #3 \advance\count107 by #5 \divide\count107 by 3 
\advance\count108 by #4 \advance\count108 by #6 \divide\count108 by 3
\LiNe{#1}{#2}{\count107}{\count108}{#7} 
\LiNe{#3}{#4}{\count107}{\count108}{#7} 
\LiNe{#5}{#6}{\count107}{\count108}{#7} 
}}

\def\trigg#1#2#3#4#5#6#7{{\LiNe{#1}{#2}{#3}{#4}{#7} 
\LiNe{#3}{#4}{#5}{#6}{#7}\LiNe{#5}{#6}{#1}{#2}{#7}}}

\def\power#1#2{\count91=1 \count92=#2 \loop 
\advance\count92 by-1 \multiply\count91 by #1 
\ifnum \count92>0 \repeat}


\def\pciSigHanother#1#2{ \power{6}{#1}
\loop 
\advance\count91 by-1 
\count92=#1 \count95=\count91 
\count211=0
\count212=80000
\count213=160000
\count221=0
\count222=150000
\count223=0
{\loop 
\advance\count92 by-1 
\count93=\count95 
\divide\count95 by 6 
\count94=\count95 
\multiply\count94 by 6 
\advance\count93 by-\count94 
\pciSigooHanother{\count93} 
\ifnum \count92>0 
\repeat \trigg{\count211}{\count221}{\count212}{\count222}{\count213}{\count223}{#2}}
\ifnum \count91>0 
\repeat} 



\def\pciSigooHonemore#1{
\ifcase#1 
%
%
%
\count141=0\count142=0\count143=0\count151=0\count152=0\count153=0
\count231=12 
\count232=0 
\count233=0 
\multiply\count231 by \count211 
\multiply\count232 by \count212 
\multiply\count233 by \count213 
\advance\count141 by \count231
\advance\count141 by \count232
\advance\count141 by \count233
\divide\count141 by 12 
\count231=12 
\count232=0 
\count233=0 
\multiply\count231 by \count221 
\multiply\count232 by \count222 
\multiply\count233 by \count223 
\advance\count151 by \count231
\advance\count151 by \count232
\advance\count151 by \count233
\divide\count151 by 12
\count231=4 
\count232=4 
\count233=4 
\multiply\count231 by \count211 
\multiply\count232 by \count212 
\multiply\count233 by \count213 
\advance\count142 by \count231
\advance\count142 by \count232
\advance\count142 by \count233
\divide\count142 by 12 
\count231=4 
\count232=4 
\count233=4 
\multiply\count231 by \count221 
\multiply\count232 by \count222 
\multiply\count233 by \count223 
\advance\count152 by \count231
\advance\count152 by \count232
\advance\count152 by \count233
\divide\count152 by 12
\count231=5 
\count232=5 
\count233=2 
\multiply\count231 by \count211 
\multiply\count232 by \count212 
\multiply\count233 by \count213 
\advance\count143 by \count231
\advance\count143 by \count232
\advance\count143 by \count233
\divide\count143 by 12 
\count231=5 
\count232=5 
\count233=2 
\multiply\count231 by \count221 
\multiply\count232 by \count222 
\multiply\count233 by \count223 
\advance\count153 by \count231
\advance\count153 by \count232
\advance\count153 by \count233
\divide\count153 by 12
\or
\count141=0\count142=0\count143=0\count151=0\count152=0\count153=0
\count231=12 
\count232=0 
\count233=0 
\multiply\count231 by \count211 
\multiply\count232 by \count212 
\multiply\count233 by \count213 
\advance\count141 by \count231
\advance\count141 by \count232
\advance\count141 by \count233
\divide\count141 by 12 
\count231=12 
\count232=0 
\count233=0 
\multiply\count231 by \count221 
\multiply\count232 by \count222 
\multiply\count233 by \count223 
\advance\count151 by \count231
\advance\count151 by \count232
\advance\count151 by \count233
\divide\count151 by 12
\count231=4 
\count232=4 
\count233=4 
\multiply\count231 by \count211 
\multiply\count232 by \count212 
\multiply\count233 by \count213 
\advance\count142 by \count231
\advance\count142 by \count232
\advance\count142 by \count233
\divide\count142 by 12 
\count231=4 
\count232=4 
\count233=4 
\multiply\count231 by \count221 
\multiply\count232 by \count222 
\multiply\count233 by \count223 
\advance\count152 by \count231
\advance\count152 by \count232
\advance\count152 by \count233
\divide\count152 by 12
\count231=5 
\count232=2 
\count233=5 
\multiply\count231 by \count211 
\multiply\count232 by \count212 
\multiply\count233 by \count213 
\advance\count143 by \count231
\advance\count143 by \count232
\advance\count143 by \count233
\divide\count143 by 12 
\count231=5 
\count232=2 
\count233=5 
\multiply\count231 by \count221 
\multiply\count232 by \count222 
\multiply\count233 by \count223 
\advance\count153 by \count231
\advance\count153 by \count232
\advance\count153 by \count233
\divide\count153 by 12
\or
\count141=0\count142=0\count143=0\count151=0\count152=0\count153=0
\count231=2 
\count232=5 
\count233=5 
\multiply\count231 by \count211 
\multiply\count232 by \count212 
\multiply\count233 by \count213 
\advance\count141 by \count231
\advance\count141 by \count232
\advance\count141 by \count233
\divide\count141 by 12 
\count231=2 
\count232=5 
\count233=5 
\multiply\count231 by \count221 
\multiply\count232 by \count222 
\multiply\count233 by \count223 
\advance\count151 by \count231
\advance\count151 by \count232
\advance\count151 by \count233
\divide\count151 by 12
\count231=0
\count232=12
\count233=0
\multiply\count231 by \count211 
\multiply\count232 by \count212 
\multiply\count233 by \count213 
\advance\count142 by \count231
\advance\count142 by \count232
\advance\count142 by \count233
\divide\count142 by 12 
\count231=0
\count232=12
\count233=0
\multiply\count231 by \count221 
\multiply\count232 by \count222 
\multiply\count233 by \count223 
\advance\count152 by \count231
\advance\count152 by \count232
\advance\count152 by \count233
\divide\count152 by 12
\count231=4 
\count232=4 
\count233=4 
\multiply\count231 by \count211 
\multiply\count232 by \count212 
\multiply\count233 by \count213 
\advance\count143 by \count231
\advance\count143 by \count232
\advance\count143 by \count233
\divide\count143 by 12 
\count231=4 
\count232=4 
\count233=4 
\multiply\count231 by \count221 
\multiply\count232 by \count222 
\multiply\count233 by \count223 
\advance\count153 by \count231
\advance\count153 by \count232
\advance\count153 by \count233
\divide\count153 by 12
\or
\count141=0\count142=0\count143=0\count151=0\count152=0\count153=0
\count231=5 
\count232=5 
\count233=2 
\multiply\count231 by \count211 
\multiply\count232 by \count212 
\multiply\count233 by \count213 
\advance\count141 by \count231
\advance\count141 by \count232
\advance\count141 by \count233
\divide\count141 by 12 
\count231=5 
\count232=5 
\count233=2 
\multiply\count231 by \count221 
\multiply\count232 by \count222 
\multiply\count233 by \count223 
\advance\count151 by \count231
\advance\count151 by \count232
\advance\count151 by \count233
\divide\count151 by 12
\count231=0
\count232=12
\count233=0
\multiply\count231 by \count211 
\multiply\count232 by \count212 
\multiply\count233 by \count213 
\advance\count142 by \count231
\advance\count142 by \count232
\advance\count142 by \count233
\divide\count142 by 12 
\count231=0
\count232=12
\count233=0
\multiply\count231 by \count221 
\multiply\count232 by \count222 
\multiply\count233 by \count223 
\advance\count152 by \count231
\advance\count152 by \count232
\advance\count152 by \count233
\divide\count152 by 12
\count231=4 
\count232=4 
\count233=4 
\multiply\count231 by \count211 
\multiply\count232 by \count212 
\multiply\count233 by \count213 
\advance\count143 by \count231
\advance\count143 by \count232
\advance\count143 by \count233
\divide\count143 by 12 
\count231=4 
\count232=4 
\count233=4 
\multiply\count231 by \count221 
\multiply\count232 by \count222 
\multiply\count233 by \count223 
\advance\count153 by \count231
\advance\count153 by \count232
\advance\count153 by \count233
\divide\count153 by 12
\or
\count141=0\count142=0\count143=0\count151=0\count152=0\count153=0
\count231=4
\count232=4
\count233=4
\multiply\count231 by \count211 
\multiply\count232 by \count212 
\multiply\count233 by \count213 
\advance\count141 by \count231
\advance\count141 by \count232
\advance\count141 by \count233
\divide\count141 by 12 
\count231=4 
\count232=4 
\count233=4 
\multiply\count231 by \count221 
\multiply\count232 by \count222 
\multiply\count233 by \count223 
\advance\count151 by \count231
\advance\count151 by \count232
\advance\count151 by \count233
\divide\count151 by 12
\count231=2
\count232=5
\count233=5
\multiply\count231 by \count211 
\multiply\count232 by \count212 
\multiply\count233 by \count213 
\advance\count142 by \count231
\advance\count142 by \count232
\advance\count142 by \count233
\divide\count142 by 12 
\count231=2
\count232=5
\count233=5
\multiply\count231 by \count221 
\multiply\count232 by \count222 
\multiply\count233 by \count223 
\advance\count152 by \count231
\advance\count152 by \count232
\advance\count152 by \count233
\divide\count152 by 12
\count231=0 
\count232=0 
\count233=12 
\multiply\count231 by \count211 
\multiply\count232 by \count212 
\multiply\count233 by \count213 
\advance\count143 by \count231
\advance\count143 by \count232
\advance\count143 by \count233
\divide\count143 by 12 
\count231=0 
\count232=0 
\count233=12 
\multiply\count231 by \count221 
\multiply\count232 by \count222 
\multiply\count233 by \count223 
\advance\count153 by \count231
\advance\count153 by \count232
\advance\count153 by \count233
\divide\count153 by 12
\or
\count141=0\count142=0\count143=0\count151=0\count152=0\count153=0
\count231=4
\count232=4
\count233=4
\multiply\count231 by \count211 
\multiply\count232 by \count212 
\multiply\count233 by \count213 
\advance\count141 by \count231
\advance\count141 by \count232
\advance\count141 by \count233
\divide\count141 by 12 
\count231=4 
\count232=4 
\count233=4 
\multiply\count231 by \count221 
\multiply\count232 by \count222 
\multiply\count233 by \count223 
\advance\count151 by \count231
\advance\count151 by \count232
\advance\count151 by \count233
\divide\count151 by 12
\count231=5
\count232=2
\count233=5
\multiply\count231 by \count211 
\multiply\count232 by \count212 
\multiply\count233 by \count213 
\advance\count142 by \count231
\advance\count142 by \count232
\advance\count142 by \count233
\divide\count142 by 12 
\count231=5
\count232=2
\count233=5
\multiply\count231 by \count221 
\multiply\count232 by \count222 
\multiply\count233 by \count223 
\advance\count152 by \count231
\advance\count152 by \count232
\advance\count152 by \count233
\divide\count152 by 12
\count231=0 
\count232=0 
\count233=12 
\multiply\count231 by \count211 
\multiply\count232 by \count212 
\multiply\count233 by \count213 
\advance\count143 by \count231
\advance\count143 by \count232
\advance\count143 by \count233
\divide\count143 by 12 
\count231=0 
\count232=0 
\count233=12 
\multiply\count231 by \count221 
\multiply\count232 by \count222 
\multiply\count233 by \count223 
\advance\count153 by \count231
\advance\count153 by \count232
\advance\count153 by \count233
\divide\count153 by 12
\fi%
\count211=\count141
\count212=\count142
\count213=\count143
\count221=\count151
\count222=\count152
\count223=\count153
}

\def\trigg#1#2#3#4#5#6#7{{
\count107=#1 
\count108=#2 
\advance\count107 by #3 \advance\count107 by #5 \divide\count107 by 3 
\advance\count108 by #4 \advance\count108 by #6 \divide\count108 by 3
\LiNe{#1}{#2}{\count107}{\count108}{#7} 
\LiNe{#3}{#4}{\count107}{\count108}{#7} 
\LiNe{#5}{#6}{\count107}{\count108}{#7} 
}}

\def\trigg#1#2#3#4#5#6#7{{\LiNe{#1}{#2}{#3}{#4}{#7} 
\LiNe{#3}{#4}{#5}{#6}{#7}\LiNe{#5}{#6}{#1}{#2}{#7}}}

\def\power#1#2{\count91=2 \count92=#2 \loop 
\advance\count92 by-1 \multiply\count91 by #1 
\ifnum \count92>0 \repeat}


\def\pciSigHonemore#1#2{ \power{6}{#1}
\loop 
\advance\count91 by-1 
\count92=#1 \count95=\count91 
\count211=0
\count212=80000
\count213=160000
\count221=0
\count222=150000
\count223=0
{\loop 
\advance\count92 by-1 
\count93=\count95 
\divide\count95 by 6 
\count94=\count95 
\multiply\count94 by 6 
\advance\count93 by-\count94 
\pciSigooHonemore{\count93} 
\ifnum \count92>0 
\repeat \trigg{\count211}{\count221}{\count212}{\count222}{\count213}{\count223}{#2}}
\ifnum \count91>0 
\repeat} 

\BFIG{figpciSigHanother}{A regular  post-critically infinite fractal and its first approximation.}{\begin{picture}(160,170)(0,-10)
\hhh\end{picture}
\ 
\hfil
\ 
\begin{picture}(160,170)(0,-10)
\end{picture}}
\BFIG{figpciSigHonemore}{A  non regular  post-critically infinite fractal and its first approximation.}{\begin{picture}(160,170)(0,-10)
\hhh\end{picture}
\ 
\hfil
\ 
\begin{picture}(160,170)(0,-10)
\end{picture}}


\BEX{e-nSG}{Post-critically infinite \Sig}{
The  post-critically infinite \Sig\ is a \frf\ which has many properties of the \Sig, but is not  a \pcfsss. 
More exactly, its post-critical set defined in \cite {Ki1,Ki4} is countably infinite, and each 
vertex $v\in V_*$ is an intersection of countably
many cells with pairwise disjoint interior. This fractal satisfies \Def{definition-sspcnf} and 
can be constructed as a self-affine fractal in $\mathbb R^2$ using nine contractions, see \Fig{figpciSig}. In \Fig{figpciSig} we also sketch the first approximation to it in \hc. 
In particular, \Fig{figpciSig} shows the values of a symmetric and a skew-symmetric harmonic functions. 
By \Thm{thm-symm} one can easily construct a resistance form such that for any 
$n$ the resistances are equal to $\left(50/53\right)^n$ in each triangle 
with vertices in $V_{n}$. The energy renormalization factor is $53/50=\rho_1=...=\rho_9$. The fact that this 
factor is larger than one is significant because it implies that the harmonic structure is regular by \Thm{thm-reg}, that is $\Omega=F$. 
By \Thm{cor-HC} this resistance form 
 satisfies all the assumptions, including the (HC) assumption. }

\BEX{e-nSGtwomore}{}{In the end we describe two more examples of 
  post-critically infinite  \frf s, which are shown in Figures~\ref{figpciSigHanother} and \ref{figpciSigHonemore}. In these examples for any $n$ there are $n$-cells which are 
joined in two points.  Both fractals 
satisfy \Def{definition-sspcnf} and 
can be constructed as a self-affine fractal in $\mathbb R^2$ using six contractions. In Figures~\ref{figpciSigHanother} and \ref{figpciSigHonemore}  we also sketch the first approximations to these fractals in \hc. 
In particular, one can see the values of  symmetric and  skew-symmetric harmonic functions on each fractal. 
By \Thm{thm-symm} one can easily construct  resistance forms such that $E_0$ is given by~\eqref{e-e0}. By \Thm{cor-HC} these resistance forms 
 satisfy  the (HC) assumption. 
In the case of the fractal in Figures~\ref{figpciSigHanother}, an elementary calculation shows that the common energy renormalization factor in~\eqref{e-rho}  is $5/4$, and so the resistance form is regular. 
In the case of the fractal in Figures~\ref{figpciSigHonemore}, the calculation shows that the common energy renormalization factor in~\eqref{e-rho}  is $4/5$, and so the resistance form is non regular.}

\BREM{rem-spectrdec}{If the  assumptions of \Thm{thm-symm} are satisfied and  a  \Lp\ is defined with respect to the product (Bernoulli) measure that gives equal weight to all $n$-cells, then   one can compute the spectrum of this \Lp\ by the so called spectral decimation method  of \cite{MT,T00}. In particular, this can be done for the  fractals shown in Figures~\ref{figpciSig}, \ref{figpciSigHanother} and~\ref{figpciSigHonemore}. Note, however, that the results of \Sec{section-Lp} are not applicable to such a \Lp.}

\end{document}